\documentclass[a4paper,11pt,twoside]{article}
\usepackage[top=3cm, bottom=3cm, inner=2cm, outer=2cm, includehead]{geometry}
\usepackage{fancyhdr}
\usepackage{xurl}
\usepackage{graphicx}
\usepackage{alltt}
\usepackage[hidelinks,hypertexnames=false]{hyperref}
\usepackage{lmodern}
\usepackage{csquotes}
\usepackage{amsfonts}
\usepackage[sort,nocompress]{cite}
\newcommand{\autocite}{\cite}
\providecommand{\RN}[1]{\textup{\uppercase\expandafter{\romannumeral#1}}}
\usepackage{romannum}
\usepackage[english]{babel}
\usepackage{setspace}
\usepackage{amsmath}
\usepackage{enumitem}
\setlist[itemize]{leftmargin=0pt,itemindent=1.8em,labelsep=.5em,nosep,topsep=3pt}
\usepackage{booktabs}
\usepackage{float}
\graphicspath{{figures/}}
\usepackage{amsthm}     
\usepackage{amssymb}  
\usepackage{siunitx}   
\usepackage{booktabs}  
\usepackage{graphicx}     
\usepackage{caption}      
\graphicspath{{figures/}} 

\theoremstyle{definition} 
\newtheorem{theorem}{Theorem}[section]  

\theoremstyle{definition}  
\newtheorem{definition}[theorem]{Definition}  

\theoremstyle{definition}  
\newtheorem{assumption}[theorem]{Assumption}  

\theoremstyle{definition}  
\newtheorem{lemma}[theorem]{Lemma}  

\theoremstyle{definition}  
\newtheorem{corollary}[theorem]{Corollary}  

\theoremstyle{definition}  
\newtheorem{proposition}[theorem]{Proposition}  
\newtheorem*{remark}{Remark}
\usepackage{xcolor}
\hypersetup{
  pdftitle={Discrete-Time Zero-Sum Markov Games with Control, Stopping, and State-Dependent Discounting},
  pdfauthor={Guanghui Feng and Xin Guo}
}
\begin{document}
	\pagenumbering{arabic}
	\setcounter{page}{1}
	\onehalfspacing
	
	\begin{center}
		\LARGE
		\textbf{Discrete-Time Zero-Sum Markov Games with Control, Stopping, and State-Dependent Discounting}\\[12pt]
		\normalsize
		\textbf{Guanghui Feng}$^{1}$\qquad\textbf{Xin Guo}$^{2,*}$\\[8pt]
		\begin{singlespace}
		\small
		$^{1}$School of Mathematics, Sun Yat-sen University, China\\[4pt]
		$^{2}$School of Science, Sun Yat-sen University, China
		\end{singlespace}
	\end{center}
	\begingroup
	\renewcommand{\thefootnote}{\fnsymbol{footnote}}
	\footnotetext[1]{Corresponding author at: School of Science, Sun Yat-sen University,
	66 Gongchang Road, Guangming District, Shenzhen 518107, China.
	\par\noindent E-mail addresses:
	\href{mailto:fenggh7@mail2.sysu.edu.cn}{fenggh7@mail2.sysu.edu.cn} (G. Feng),
	\href{mailto:guox87@mail.sysu.edu.cn}{guox87@mail.sysu.edu.cn} (X. Guo).}
	\endgroup

	\begin{abstract}
		\normalsize
		We study a discrete-time two-player zero-sum stochastic game with control and stopping, where both players can stop the game, the discount factors depends on the states, and the running and terminal rewards may be unbounded. A contraction argument for the constrained Shapley equation establishes the existence of the value and a Nash equilibrium consisting of stationary controls and stopping times. We further provides algorithms for approximating the value and computing an $\varepsilon$-Nash equilibrium. A numerical example with state-dependent discount factors illustrates the results.

	\end{abstract}

	\section{Introduction}\label{s:1}
	Dynamic game problems often involve both the operation of a system and the decision to terminate it. In an economic or operational relationship, the parties may adjust their actions while the system continues, but either party may also wish to end the relationship. Controls then determine the evolution of a state in the system, whereas stopping decisions determine how long rewards accumulate. The importance of future rewards may also vary with the current economic or operating state, motivating state-dependent discounting \autocite{Wei_Guo_2011,Guo_Chen_Li_2023}. For example, if a firm's financing cost varies with its current cash reserves, the discount factor used to evaluate future returns may also depend on the current state. Moreover, rewards that grow with the size of the system need not be uniformly bounded. In a population model, for example, a total service cost proportional to the population size is unbounded when the population can take arbitrarily large integer values. These considerations lead to the problem studied here: a discrete-time zero-sum game in which both players control and stop, with state-dependent discount factors and possibly unbounded running and terminal rewards. For applications, we not only establish an equilibrium, but also give the algorithms to compute an $\varepsilon$-Nash equilibrium.

	We first review related studies on stochastic games involving control, stopping, or both. In continuous time, stopping games allow players to choose when to stop a given process without controlling its dynamics \autocite{DeAngelis_Ferrari_Moriarty_2018,Christensen_Lindensjo_2026,Perez_Rodosthenous_Yamazaki_2025}. When both players can control and stop, Pal and Saha \autocite{Pal_Saha_2020} study a continuous-time Markov chain on a countable state space with a constant discount rate and possibly unbounded transition and reward rates. They establish the value and a saddle-point equilibrium. Related semi-Markov models also allow both players to control and stop. Bose, Pal, and Saha \autocite{Bose_Pal_Saha_2026} consider a risk-sensitive discounted criterion in a semi-Markov setting and establish the value and a saddle-point equilibrium. In discrete time, classical stopping games, beginning with Dynkin \autocite{Dynkin_1969}, allow players to choose when to stop a given process without controlling its transition law. Subsequent studies establish value and equilibrium results under different reward assumptions and game structures \autocite{Rosenberg_Solan_Vieille_2001,Shmaya_Solan_2004,Hamadene_Hassani_2014}. A further class assigns control to one player and stopping to the other. Cavazos-Cadena and Hernández-Hernández \autocite{Cavazos_Hernandez_2012} study such a game under an expected discounted criterion, establish the value and equilibria, and construct approximate equilibria through successive approximations. Under an expected total-reward criterion, equilibrium results are obtained under absorption conditions \autocite{Cavazos_Rodriguez_Sanchez_2021} or communication and ergodicity conditions \autocite{Cavazos_Cantu_Cerda_2021}. More directly related to this paper are discrete-time games in which both players can control and stop. Ghosh and Mallikarjuna Rao \autocite{Ghosh_2007} consider a Borel state space, constant discounting, bounded continuous rewards, and strictly ordered terminal rewards. They characterize the value by dynamic programming inequalities and establish a saddle-point equilibrium. In contrast, the present paper considers an expected discounted criterion with state-dependent discount factors and possibly unbounded running and terminal rewards.

	State-dependent discounting has also been studied without this two-player control-and-stopping structure. For a single decision maker, Wei and Guo \autocite{Wei_Guo_2011} allow unbounded rewards or costs, establish optimal stationary policies, and develop value iteration. Jasso-Fuentes, Menaldi, and Prieto-Rumeau \autocite{Jasso_Menaldi_Prieto_2020} allow state-action-dependent discount factors, and characterize optimal policies and the value by dynamic programming; their framework includes control problems with stopping. For stochastic games without strategic stopping, Yu, Guo, and Xia \autocite{Yu_Guo_Xia_2022} establish the value and optimal stationary strategies under an expected criterion with state-action-dependent discounting and weighted reward bounds, and develop an algorithm for the value and approximate equilibria in finite models. Guo, Chen, and Li \autocite{Guo_Chen_Li_2023} address discrete-time risk-sensitive games with varying discount factors and unbounded rewards under logarithmic growth conditions, and provide existence results and an iterative algorithm. These works supply important methods for variable discounting, but do not determine both players' control policies and stopping times for the expected criterion studied here.

	The gap addressed in this paper is the joint treatment of these decisions and reward conditions, together with an equilibrium computation algorithm. The expected-discounted control-and-stopping studies cited above use constant discount factors or rates, whereas the variable-discount control results concern a single decision maker or games without both players' stopping decisions. These results do not directly give the existence and computation guarantees under the combined assumptions considered here. In discrete time, state-dependent discounting requires more than a notational replacement: the cumulative discount depends on the realized state path, whose distribution is affected by both controls. With unbounded rewards, one must control the discounted running rewards and the terminal rewards at random stopping times, and show how errors in the value approximation affect both controls and stopping decisions.

	We establish the existence of the value and a Nash equilibrium on a countable state space under explicit reward-growth and transition conditions, with state-dependent discount factors uniformly bounded above by a constant less than one. First, a contraction argument in a weighted norm gives a unique solution of the constrained Shapley equation despite possibly unbounded rewards. Second, bounds for the cumulative discounts and terminal terms allow us to identify this solution with the value and construct a Nash equilibrium consisting of stationary controls and stopping times. These results treat state-dependent discounting under an expected criterion, in contrast to the constant-discount formulation in \autocite{Ghosh_2007} and the risk-sensitive formulation in \autocite{Zhang_Liu_2024}. Third, following the contraction-based value-iteration method in Puterman \autocite[Sections 6.2.3 and 6.3.2]{Puterman_1994}, we obtain an error bound from successive terms of a single sequence. We first prove a weighted error bound and then, under the finiteness condition in Section 5, derive the uniform value-error bound needed for the Nash inequalities. The value-iteration principle is classical; its use here requires error bounds that account for both controls and stopping decisions under state-dependent discounting. For finite state and action spaces, the resulting algorithms compute the value approximation, both control policies, and both stopping sets by solving linear programs. Together, these results connect existence, value approximation, and equilibrium computation for a model combining both players' control and stopping decisions with state-dependent discounting and possibly unbounded rewards. An infinite-state queueing example illustrates the computed strategies and the error guarantees.

	The remainder of the paper is organized as follows. Section~\ref{s:2} formulates the game and its criterion. Section~\ref{s:3} establishes the existence and uniqueness of the solution to the constrained Shapley equation. Section~\ref{s:4} proves the existence of the value and a Nash equilibrium. Section~\ref{s:5} develops the approximation procedures and algorithms, and Section~\ref{s:6} presents a numerical example.
	
	\section{Problem formulation}\label{s:2}
	For a countable discrete state space $X$, let $\mathcal{B}(X)$ denote its Borel $\sigma$-algebra, and let $\mathcal{P}(X)$ denote the set of probability measures on $\mathcal{B}(X)$, endowed with the topology of weak convergence.
	
	A discrete-time two-player zero-sum stochastic game with control and stopping is
	\begin{align}
		\{X, (A(x)\subset A,B(x)\subset B,x\in X), p, r(x,a,b),g_1(x), g_2(x)\}. \tag{2.1}
	\end{align}
	
	$X$ is a countable state space; $A(x)\in \mathcal{B}(A)$ and $B(x)\in \mathcal{B}(B)$ are the control sets available to players 1 and 2, respectively, in state $x$; $A$ and $B$ are the corresponding Borel action spaces; $p$ is a stochastic kernel on $X$ given $K := \{(x,a,b) \mid x\in X,a\in A(x),b\in B(x)\}$ that describes the transition law. For each $(x,a,b)\in K$, $p(\cdot\mid x,a,b)$ is a probability measure on $X$; for each $y\in X$, $p(y\mid\cdot)$ is a real-valued Borel-measurable function on $K$; $r\colon K\to\mathbb{R}$ is a measurable one-stage reward function; and $g_k\colon X\to\mathbb{R}$ is the terminal payoff function when player $k$ stops the game, for $k=1,2$.

	For the model in (2.1), let $H_t$ be the space of admissible histories up to time $t\in N_0=\{0,1,\ldots\}$. That is, $H_0:=X$ and, for $t=1,2,\ldots$, 
	\begin{equation}
		H_t := K^t \times X = K \times H_{t-1} \ \ \ t=1,2,\ldots. \tag{2.2}
	\end{equation} 
	
	A randomized control policy for player 1 is a sequence $\pi^1=\{\pi^1_t:t\in N_0\}$ of stochastic kernels $\pi^1_t\colon H_t\to\mathcal{P}(A)$ satisfying $\pi^1_t(A(x_t)| h_t) = 1$ for all $h_t:=(x_0,a_0,b_0,\ldots,x_{t-1},a_{t-1},b_{t-1},x_t)\in H_t, t\in N_0.$ The policy $\pi^1$ is stationary if there is a stochastic kernel $\varphi^1$ on $A$ given $X$ such that $\pi^1_t(\cdot| h_t) = \varphi^1(\cdot| x_t)$ for all $h_t:=(x_0,a_0,b_0,\ldots,x_{t-1},a_{t-1},b_{t-1},x_t)\in H_t, t\in N_0.$

	Let $\Pi^1$ and $\Pi^1_s$ denote, respectively, the sets of all policies and all stationary policies for player 1. The policies for player 2 are defined analogously; with $B(x)$ in place of $A(x)$, we denote the corresponding sets by $\Pi^2$ and $\Pi^2_s$.
	
	Let $\Omega=(X\times A\times B)^\infty$, equipped with a $\sigma$-algebra $\mathcal{F}$. Define the random variables $X_t$, $A_t$, $B_t$, and the random history $Z_t$ on $(\Omega,\mathcal{F})$ by
	\begin{equation}
		X_t(\omega)=x_t, \quad A_t(\omega)=a_t, \quad B_t(\omega)=b_t, \quad
		Z_t(\omega)=(x_0,a_0,b_0,\ldots,x_{t-1},a_{t-1},b_{t-1},x_t)=:h_t, \tag{2.3}
	\end{equation}
	for $t=0,1,2,\ldots$, where $\omega = (x_0, a_0, b_0, x_1, a_1, b_1,\ldots) \in \Omega$; thus $Z_t(\omega)\in H_t$.
	
	Given an initial distribution $\nu$ on $X$ and a policy pair $(\pi^1,\pi^2)\in\Pi^1\times\Pi^2$, the Ionescu--Tulcea theorem yields a unique probability measure $P^{\pi^1,\pi^2}_\nu$ on $(\Omega,\mathcal{F})$ such that, for every $G\in\mathcal{B}(X)$, $C_1\in\mathcal{B}(A)$, $C_2\in\mathcal{B}(B)$, and $h_t\in H_t$, $P^{\pi^1,\pi^2}_\nu\{X_0\in G\}=\nu(G),$
	\begin{align}
			P^{\pi^1,\pi^2}_\nu
			\{A_t\in C_1,B_t\in C_2\mid Z_t=h_t\}
			=\pi_t^1(C_1\mid h_t)\pi_t^2(C_2\mid h_t),& \tag{2.4}\\
			P^{\pi^1,\pi^2}_\nu
			\{X_{t+1}\in G\mid Z_t=h_t,A_t=a_t,B_t=b_t\}
			=p(G\mid x_t,a_t,b_t).& \tag{2.5}
	\end{align}
	
	Let \(E^{\pi^1,\pi^2}_\nu\) denote expectation with respect to \(P^{\pi^1,\pi^2}_\nu\). When the initial distribution \(\nu\) is concentrated at \(x\), we denote the corresponding probability measure and expectation by \(P^{\pi^1,\pi^2}_x\) and \(E^{\pi^1,\pi^2}_x\), respectively. To describe stoppings, let $\mathcal{G}_t:=\sigma(Z_t), t\in N_0$ be the $\sigma$-field generated by $Z_t$. A random time $\theta\colon\Omega\to\mathbb{N}_0\cup\{+\infty\}$ is a stopping time if $\{\theta\leq t\}\in\mathcal{G}_t$ for every $t\in\mathbb{N}_0$. Let $\Theta$ denote the set of all stopping times.
	\begin{definition}
		A strategy for player $k$ is a pair $(\pi^k, \theta^k)$, where $\pi^k \in\Pi^k$ and $\theta^k \in\Theta$.
	\end{definition}
	Let $D\subset X$ be a prescribed set of states at which the system fails. Define the time of system failure by $\sigma_D:=\inf\{t\in\mathbb N_0:X_t\in D\},$ where $\inf\varnothing:=\infty.$
	Since
	\[
		\{\sigma_D\leq t\}=\bigcup_{s=0}^{t}\{X_s\in D\}\in\mathcal G_t,
		\qquad t\in\mathbb N_0,
	\]
	$\sigma_D$ is a stopping time with respect to $(\mathcal G_t)$. Unlike the stopping times chosen by the players, $\sigma_D$ is determined by the state process and the prescribed set $D$.

	The game proceeds as follows. Given stopping times $\theta^1$ and $\theta^2$ chosen by the two players, at each time $t$, the players observe the current state $x_t$ and the previously realized actions. If the current state belongs to $D$, the game terminates immediately, with no additional terminal reward. Otherwise, the players independently decide whether to stop or continue. If both players continue, players 1 and 2 select $a_t\in A(x_t)$ and $b_t\in B(x_t)$, respectively. Player 1 then receives the one-stage payoff $r(x_t,a_t,b_t)$ from player 2, and the state moves to $x_{t+1}\in X$ according to $p(x_{t+1}\mid x_t,a_t,b_t)$. If player 1 stops at time $t$ before system failure and no later than player 2, player 1 receives the terminal reward $g_1(x_t)$.
	If player 2 stops at time $t$ before system failure and before
	player 1, player 1 receives the terminal reward $g_2(x_t)$. Thus all rewards in the criterion below are measured from player 1's perspective. Player 1 seeks to maximize, and player 2 to minimize, the expected discounted criterion. The game terminates at $\theta^1\wedge\theta^2\wedge\sigma_D$. If system failure and a player's stopping decision occur at the same time, the game terminates according to the system-failure rule. This procedure is repeated at each subsequent time until either player stops or the system fails.
	
	For each state $x\in X$, let $\alpha(x)\in(0,1)$ be the corresponding discount factor. Define the discounted criterion, provided it is well defined, by
	\begin{align}
		J(x,& \pi^1, \theta^1, \pi^2, \theta^2)
		=E_x^{\pi^1,\pi^2}\bigg[
		\sum_{t=0}^{(\theta^1\wedge\theta^2\wedge\sigma_D)-1}
		\left(\prod_{k=0}^{t-1}\alpha(X_k)\right)r(X_t,A_t,B_t) \notag\\
		&+\left(\prod_{k=0}^{\theta^1-1}\alpha(X_k)\right)g_1(X_{\theta^1})
		I_{\{\theta^1\leq\theta^2,\,\theta^1<\sigma_D\}}
		+\left(\prod_{k=0}^{\theta^2-1}\alpha(X_k)\right)g_2(X_{\theta^2})
		I_{\{\theta^2<\theta^1,\,\theta^2<\sigma_D\}}
		\bigg]. \tag{2.6}
	\end{align}
	\begin{remark}
		If $\alpha(x)\equiv\alpha_0\in(0,1)$, then (2.6) reduces to the constant-discount criterion, with cumulative discount $\alpha_0^t$, as in the constant-discount formulation of \autocite{Ghosh_2007}. If $D=\varnothing$, then $\sigma_D=\infty$ and there is no system-failure termination, as in the game studied by \autocite{Wu_Wang_Kong_2021}.
	\end{remark}
	For $x\in D$, $J(x, \pi^1, \theta^1, \pi^2, \theta^2)=0$ for every strategy pair. Henceforth, we consider initial states in $D^c:=X\setminus D$; functions representing the value or its approximations are extended by zero to $D$ whenever evaluated along the state process. The functions on $D^c$ defined as
	\begin{align}
		J^-(x) := \sup_{(\pi^1, \theta^1)\in \Pi^1\times \Theta}\inf_{(\pi^2, \theta^2)\in \Pi^2\times \Theta} J(x, \pi^1, \theta^1, \pi^2, \theta^2) \tag{2.7}, \\
		J^+(x) := \inf_{(\pi^2, \theta^2)\in \Pi^2\times \Theta}\sup_{(\pi^1, \theta^1)\in \Pi^1\times \Theta} J(x, \pi^1, \theta^1, \pi^2, \theta^2) \tag{2.8}
	\end{align}
	are called the lower and upper values of the game, respectively. Clearly, $J^-(x)\leq J^+(x)$.
	\begin{definition}
		If $J^-(x)=J^+(x)$ for every $x\in D^c$, then their common function, denoted by $J$, is called the value of the game.
	\end{definition}
	\begin{definition}
		Suppose that the game has a value $J$. A strategy $(\pi^1_*,\theta^1_*)$ is optimal for player 1 if 
		\begin{equation}
			\inf_{(\pi^2, \theta^2)\in \Pi^2\times \Theta}J(x, \pi^1_*, \theta^1_*, \pi^2, \theta^2) = J(x) \ \ \ \forall x\in D^c. \tag{2.9}
		\end{equation}
		Similarly, $(\pi^2_*,\theta^2_*)$ is optimal for player 2 if
		\begin{equation}
			\sup_{(\pi^1, \theta^1)\in \Pi^1\times \Theta} J(x, \pi^1, \theta^1, \pi^2_*, \theta^2_*) = J(x) \ \ \ \forall x\in D^c. \tag{2.10}
		\end{equation}
		If both players have optimal strategies $(\pi^1_*,\theta^1_*)$ and $(\pi^2_*,\theta^2_*)$, then $(\pi^1_*,\theta^1_*;\pi^2_*,\theta^2_*)$ is a Nash equilibrium.
		
		For $\varepsilon>0$, a strategy pair $(\pi^1_\varepsilon,\theta^1_\varepsilon;\pi^2_\varepsilon,\theta^2_\varepsilon)$ is an $\varepsilon$-Nash equilibrium if 
		\begin{equation}
			\inf_{(\pi^2, \theta^2)\in \Pi^2\times \Theta}J(x, \pi^1_\varepsilon, \theta^1_\varepsilon, \pi^2, \theta^2) \geq J(x)-\varepsilon \ \ \ \forall x\in D^c, \tag{2.11}
		\end{equation}
		and
		\begin{equation}
			\sup_{(\pi^1, \theta^1)\in \Pi^1\times \Theta} J(x, \pi^1, \theta^1, \pi^2_\varepsilon, \theta^2_\varepsilon) \leq J(x)+\varepsilon \ \ \ \forall x\in D^c. \tag{2.12}
		\end{equation}
	\end{definition}
	
	\section{Preliminaries}\label{s:3}
	This section develops the weighted-space framework and proves existence and uniqueness for the constrained Shapley equation. The following drift and growth condition ensures that the discounted criterion is well defined, that is, $J(x, \pi^1, \theta^1, \pi^2, \theta^2)$ is finite.
	\begin{assumption}
		There exist a real-valued function $V_0\geq1$ on $X$ and constants $c_0\geq0$ and $1\leq m<\frac{1}{\bar\alpha}$, where $\bar\alpha:=\sup_{x\in X}\alpha(x)$ and $\bar\alpha<1$ such that
		\begin{itemize}
		\item[(a)] $\displaystyle\sum_{y\in X}V_0(y)p(y\mid x,a,b)\leq mV_0(x)$ for all $(x,a,b)\in K$.\hfill\textup{(3.1)}
		\item[(b)] $\max\{|r(x,a,b)|,|g_1(x)|,|g_2(x)|\}\leq c_0V_0(x)$ for all $(x,a,b)\in K$.\hfill\textup{(3.2)}
		\end{itemize}
	\end{assumption}
	\begin{remark}
		Assumption 3.1 extends the growth conditions for Markov decision processes with unbounded rewards in \autocite{Wessels_1977,VanNunen_Wessels_1978,Puterman_1994,Wei_Guo_2011} to the present two-player control-and-stopping game. Assumption 3.1 also extends the bounded-reward conditions for Markov games in \autocite{Ghosh_2007} to the case of unbounded rewards.
	\end{remark}
	Under Assumption 3.1, we next show that $J(x, \pi^1, \theta^1, \pi^2, \theta^2)$ is absolutely integrable.
	\begin{lemma}
		For every $(\pi^1,\pi^2)\in\Pi^1\times\Pi^2$, $(\theta^1,\theta^2)\in\Theta\times\Theta$, and $x\in X$,
		\begin{align*}
		E_x^{\pi^1,\pi^2}\bigg[
		\sum_{t=0}^{(\theta^1\wedge\theta^2\wedge\sigma_D)-1}
		&\left(\prod_{k=0}^{t-1}\alpha(X_k)\right)|r(X_t,A_t,B_t)|
		+\left(\prod_{k=0}^{\theta^1-1}\alpha(X_k)\right)|g_1(X_{\theta^1})|
		I_{\{\theta^1\leq\theta^2,\,\theta^1<\sigma_D\}}\notag\\
		&+\left(\prod_{k=0}^{\theta^2-1}\alpha(X_k)\right)|g_2(X_{\theta^2})|
		I_{\{\theta^2<\theta^1,\,\theta^2<\sigma_D\}}
		\bigg]
		\leq \frac{2c_0}{1-\bar\alpha m}V_0(x).
		\end{align*}
	\end{lemma}
	\begin{proof}
		By (2.5) and Assumption 3.1(a), for every $t\geq1$,
		\begin{align*}
			E^{\pi^1,\pi^2}_x\!\left[
			V_0(X_t)\,\middle|\,
			\sigma(Z_{t-1},A_{t-1},B_{t-1})
			\right]=\sum_{y\in X}V_0(y)
			p(y\mid X_{t-1},A_{t-1},B_{t-1})
			\leq mV_0(X_{t-1}).
		\end{align*}
		Consequently,
		\begin{align*}
			E_x^{\pi^1,\pi^2}[V_0(X_t)]
			\leq mE_x^{\pi^1,\pi^2}[V_0(X_{t-1})]
			\leq m^2E_x^{\pi^1,\pi^2}[V_0(X_{t-2})]
			\leq\cdots\leq m^tV_0(x).
		\end{align*}
		Combining the Assumption 3.1(b) gives
		\begin{align*}
			E_x^{\pi^1,\pi^2}&\bigg[
			\sum_{t=0}^{(\theta^1\wedge\theta^2\wedge\sigma_D)-1}
			\left(\prod_{k=0}^{t-1}\alpha(X_k)\right)|r(X_t,A_t,B_t)|
			+\left(\prod_{k=0}^{\theta^1-1}\alpha(X_k)\right)|g_1(X_{\theta^1})|
			I_{\{\theta^1\leq\theta^2,\,\theta^1<\sigma_D\}}\\
			&+\left(\prod_{k=0}^{\theta^2-1}\alpha(X_k)\right)|g_2(X_{\theta^2})|
			I_{\{\theta^2<\theta^1,\,\theta^2<\sigma_D\}}
			\bigg]
			\leq\frac{c_0}{1-\bar\alpha m}V_0(x)
			+\frac{c_0}{1-\bar\alpha m}V_0(x)
			<\infty.\qedhere
		\end{align*}
	\end{proof}
	
	A function $u$ on $X$ is $V_0$-bounded if its weighted norm $\lVert u\rVert_{V_0}:=\sup_{x\in X}\frac{|u(x)|}{V_0(x)}$ is finite. Denote the Banach space of all such functions by $M_{V_0}(X)$.
	
	To analyze the constrained Shapley equation (3.6) below, we impose the following standard compactness and continuity condition.
	\begin{assumption}
		For every $x, y\in X$,
		\begin{itemize}
		\item[(a)] $A(x)$ and $B(x)$ are compact; 
		\item[(b)] $r(x, a, b), p(y| x, a, b)$ and $\sum_{y\in X}V_0(y)p(y| x,a,b)$ are continuous in $(a, b) \in A(x)\times B(x)$.
		\end{itemize}
	\end{assumption}
	
	With $c_0$, $m$, and $V_0$ as in Assumption 3.1, we define
	\begin{equation}
		\mathbb{I}_\alpha u(x) := \sup_{\lambda\in \mathcal{P}(A(x))}\inf_{\mu\in \mathcal{P}(B(x))}\int_{A(x)}\int_{B(x)} \mathbb{I}_\alpha u(x,a,b)\mu(db)\lambda(da) \tag{3.3}
	\end{equation}
	where $\mathbb{I}_\alpha u(x,a,b):= r(x,a,b) + \alpha(x)\sum_{y\in X} u(y)p(y|x,a,b)$
	for each $x\in D^c$ and $u\in M_{V_0}(X)$. Define the operator $\mathbb{T}_\alpha$ on $M_{V_0}(X)$ by the following formula on $D^c$, with its output extended by zero to $D$:
	\begin{equation}
		\mathbb{T}_\alpha u(x) :=\max\{\min\{\mathbb{I}_\alpha u(x),g_2(x)\},g_1(x)\},\qquad x\in D^c.\tag{3.4}
	\end{equation}
	We construct the sequence $\{\hat{u}_n\}$ recursively by
	\begin{equation}
		\hat u_0(x):=g_1(x),\quad x\in D^c,
		\qquad \hat u_n:=\mathbb T_\alpha\hat u_{n-1},\quad n\geq1.\tag{3.5}
	\end{equation}
	\begin{theorem}
		Under Assumptions 3.1 and 3.3, the following statements hold.
		\begin{itemize}
		\item[(a)] $\hat{u}_{n+1}(x)\geq \hat{u}_n(x)\geq g_1(x)$ for $x\in D^c$ and $n\geq0$. The sequence $\{\hat{u}_n\}$ converges to $u^*$ in $M_{V_0}(X)$; 
		\item[(b)] The function $u^*$ is the unique fixed point of $\mathbb{T}_\alpha$ in $M_{V_0}(X)$ and satisfies
		\[
		\left\{\begin{aligned}
			u^*(x)&=\max\{\min\{\mathbb I_\alpha u^*(x),g_2(x)\},g_1(x)\},\\
			g_1(x)&\leq u^*(x)\leq\max\{g_1(x),g_2(x)\},\qquad x\in D^c.
		\end{aligned}\right.\tag{3.6}
		\]
		\end{itemize}
	\end{theorem}
	\begin{proof}
		Set $\kappa:=\bar\alpha m<1$. We first verify that $\mathbb T_\alpha$ maps $M_{V_0}(X)$ into itself, so that its fixed point can be studied in this space. For every $u\in M_{V_0}(X)$ and $(x,a,b)\in K$ with $x\in D^c$, Assumption 3.1 gives
		\begin{align*}
			&|\mathbb I_\alpha u(x,a,b)|
			\leq |r(x,a,b)|
			+\alpha(x)\sum_{y\in X}|u(y)|p(y\mid x,a,b)\\
			&\leq c_0V_0(x)
			+\alpha(x)\lVert u\rVert_{V_0}
			\sum_{y\in X}V_0(y)p(y\mid x,a,b)
			\leq\bigl(c_0+\kappa\lVert u\rVert_{V_0}\bigr)V_0(x).
		\end{align*}
		It follows that
		\[
			|\mathbb I_\alpha u(x)|
			\leq\bigl(c_0+\kappa\lVert u\rVert_{V_0}\bigr)V_0(x).
		\]
		For $x\in D^c$, since $|g_k(x)|\leq c_0V_0(x)$ for $k=1,2$,
		\begin{align*}
			|\mathbb T_\alpha u(x)|
			\leq\max\{|\mathbb I_\alpha u(x)|,|g_1(x)|,|g_2(x)|\}
			\leq\bigl(c_0+\kappa\lVert u\rVert_{V_0}\bigr)V_0(x).
		\end{align*}
		Therefore, $\mathbb T_\alpha u\in M_{V_0}(X)$.

		We next establish the inequality needed to prove the contraction property of $\mathbb T_\alpha$. For $u,v\in M_{V_0}(X)$ and $(x,a,b)\in K$ with $x\in D^c$,
		\begin{align*}
			|\mathbb I_\alpha u(x,a,b)-\mathbb I_\alpha v(x,a,b)|
			\leq\alpha(x)\sum_{y\in X}|u(y)-v(y)|p(y\mid x,a,b)
			\leq\kappa V_0(x)\lVert u-v\rVert_{V_0}.
		\end{align*}
		Consequently,
		\begin{align*}
			|\mathbb I_\alpha u(x)-\mathbb I_\alpha v(x)|
			&\leq\sup_{\lambda\in\mathcal P(A(x))}
			\sup_{\mu\in\mathcal P(B(x))}
			\bigg|\int_{A(x)}\int_{B(x)}
			\bigl(\mathbb I_\alpha u(x,a,b)-\mathbb I_\alpha v(x,a,b)\bigr)
			\mu(db)\lambda(da)\bigg|\\
			&\leq\kappa V_0(x)\lVert u-v\rVert_{V_0}.
		\end{align*}
		For every $x\in D^c$, define
		\[
			\Gamma_x(z):=\max\{\min\{z,g_2(x)\},g_1(x)\},
			\qquad z\in\mathbb R.
		\]
		A direct calculation gives, for every $c,z,z'\in\mathbb R$,
		\[
			|\min\{z,c\}-\min\{z',c\}|\leq|z-z'|,
			\qquad
			|\max\{z,c\}-\max\{z',c\}|\leq|z-z'|.
		\]
		Hence,
		\begin{equation}
			|\Gamma_x(z)-\Gamma_x(z')| \leq |\min\{z,g_2(x)\}-\min\{z',g_2(x)\}|
			\leq|z-z'|,
			\qquad x\in D^c,\quad z,z'\in\mathbb R. \tag{3.7}
		\end{equation}
		Because $\mathbb T_\alpha u(x)=\Gamma_x(\mathbb I_\alpha u(x))$, (3.7) gives
		\[
			|\mathbb T_\alpha u(x)-\mathbb T_\alpha v(x)|
			\leq|\mathbb I_\alpha u(x)-\mathbb I_\alpha v(x)|
			\leq\kappa V_0(x)\lVert u-v\rVert_{V_0}.
		\]
		Dividing by $V_0(x)$ and using the zero extension of the operator outputs to $D$, we obtain
		\begin{equation}
			\lVert\mathbb T_\alpha u-\mathbb T_\alpha v\rVert_{V_0}
			\leq \kappa\lVert u-v\rVert_{V_0},
			\qquad u,v\in M_{V_0}(X). \tag{3.8}
		\end{equation}
		Since $\kappa<1$, $\mathbb T_\alpha$ is a contraction on $M_{V_0}(X)$.

		Having established (3.8), we now obtain the existence and uniqueness in part (b), together with the convergence in part (a). Since $M_{V_0}(X)$ is a Banach space, $\mathbb T_\alpha$ has a unique fixed point $u^*\in M_{V_0}(X)$, and, for every $v\in M_{V_0}(X)$, the sequence $\{\mathbb T_\alpha^n v\}$ converges to $u^*$ in $M_{V_0}(X)$. In particular, the sequence $\{\hat u_n\}$, where $\hat u_n=\mathbb T_\alpha^n\hat u_0$, converges to $u^*$ in $M_{V_0}(X)$.
		The fixed-point relation gives, for $x\in D^c$,
		\[
			u^*(x)
			=\max\{\min\{\mathbb I_\alpha u^*(x),g_2(x)\},g_1(x)\}.
		\]

		It remains to verify the monotonicity of the sequence in (3.5) and the second inequality in (3.6). If $u\leq v$, then $\mathbb I_\alpha u(x,a,b)\leq\mathbb I_\alpha v(x,a,b)$ because $\alpha(x)>0$ and $p$ is nonnegative. Integration, the infimum and the supremum preserve this inequality. The minimum and maximum in (3.4) preserve it on $D^c$. Thus both $\mathbb I_\alpha$ and $\mathbb T_\alpha$ are monotone. Moreover, for $x\in D^c$,
		\[
			\hat u_1(x)=\mathbb T_\alpha\hat u_0(x)\geq g_1(x)=\hat u_0(x).
		\]
		Thus,
		\[
			\hat u_2=\mathbb T_\alpha\hat u_1
			\geq\mathbb T_\alpha\hat u_0=\hat u_1,
		\]
		and induction gives
		\[
			\hat u_{n+1}(x)\geq\hat u_n(x)\geq g_1(x),
			\qquad n\geq0,\quad x\in D^c.
		\]
		For $x\in D^c$,
		\begin{align*}
			u^*(x)
			=\max\{\min\{\mathbb I_\alpha u^*(x),g_2(x)\},g_1(x)\}
			\leq\max\{g_2(x),g_1(x)\},
		\end{align*}
		because $\min\{\mathbb I_\alpha u^*(x),g_2(x)\}\leq g_2(x)$. The same fixed-point equation gives $u^*(x)\geq g_1(x)$. Therefore,
		\[
			g_1(x)\leq u^*(x)\leq\max\{g_1(x),g_2(x)\},
			\qquad x\in D^c.\qedhere
		\]
	\end{proof}
	
	\begin{remark}
		State-dependent discounting with unbounded rewards has been studied for Markov decision processes in \autocite{Wei_Guo_2011} and for semi-Markov games without stopping decisions in \autocite{Yu_Guo_Xia_2022}. The new feature relative to these studies is the treatment of state-dependent discounting in a discrete-time two-player game in which both players can control and stop.
	\end{remark}

	\section{The existence of the value and equilibrium}\label{s:4}
	This section proves the existence of the value and a Nash equilibrium. For $(x,a,b)\in K$, $\lambda\in\mathcal{P}(A(x))$, $\mu\in\mathcal{P}(B(x))$, and stationary policies $\pi^k\in\Pi_s^k$ $(k=1,2)$, define
	\begin{align}
		&p(\cdot| x, \lambda, \mu) := \int_{A(x)}\int_{B(x)}p(\cdot| x, a, b)\mu(db)\lambda(da), \ \ r(x, \lambda, \mu) := \int_{A(x)}\int_{B(x)}r(x, a, b)\mu(db)\lambda(da), \notag \\ 
		&p(\cdot| x, \pi^1, \mu) := p(\cdot| x, \varphi^1(\cdot| x), \mu), \ \ \ \ r(x, \pi^1, \mu) := r(x, \varphi^1(\cdot| x), \mu), \notag \\
		&p(\cdot| x, \lambda, \pi^2) := p(\cdot| x, \lambda, \varphi^2(\cdot| x)), \ \ \ \ r(x, \lambda, \pi^2) := r(x, \lambda, \varphi^2(\cdot| x)), \notag
	\end{align}
	where $\varphi^k(\cdot\mid x)$ is the stationary control kernel associated with $\pi^k$, $k=1,2$. Under Assumption 3.3, the functions $p(y\mid x,\lambda,\mu)$, $r(x,\lambda,\mu)$, and $\sum_{y\in X}V_0(y)p(y\mid x,\lambda,\mu)$ are continuous in $(\lambda,\mu)\in\mathcal{P}(A(x))\times\mathcal{P}(B(x))$.
	
	From the notations above, together with (3.3)
	\begin{equation}
		\mathbb{I}_\alpha u(x) := \sup_{\lambda\in \mathcal{P}(A(x))}\inf_{\mu\in \mathcal{P}(B(x))}\bigg[r(x,\lambda,\mu) + \alpha(x)\sum_{y\in X} u(y)p(y|x,\lambda,\mu)\bigg]  \ \ \forall x\in D^c.\tag{4.1}
	\end{equation}
	
	We first derive three inequalities from (3.6) which are needed to prove Theorem 4.2 below.
	\begin{lemma}
		Let $u^*$ be the solution of (3.6) in Theorem 3.4. Then the following statements hold.
		\begin{itemize}
		\item[(a)] For every $x\in D^c$ such that $g_1(x)<u^*(x)<g_2(x)$, one has $u^*(x)=\mathbb{I}_\alpha u^*(x)$; 
		\item[(b)] For every $x\in D^c$ such that $g_1(x)\leq u^*(x)<g_2(x)$, one has $\mathbb{I}_\alpha u^*(x)\leq u^*(x)$; 
		\item[(c)] For every $x\in D^c$ such that $g_1(x)<u^*(x)\leq g_2(x)$, one has $u^*(x)\leq\mathbb{I}_\alpha u^*(x)$.
		\end{itemize}
	\end{lemma}
	\begin{proof}
		Throughout this proof, $x\in D^c$.
		\begin{itemize}
		\item[(a)] Suppose that $g_1(x)<u^*(x)<g_2(x)$. Since $u^*(x)>g_1(x)$, (3.6) implies $u^*(x)=\min\{\mathbb I_\alpha u^*(x),g_2(x)\}.$ Moreover, $u^*(x)<g_2(x)$. Hence, $u^*(x)=\mathbb I_\alpha u^*(x).$

		\item[(b)] Suppose that $g_1(x)\leq u^*(x)<g_2(x)$. If $\mathbb I_\alpha u^*(x)>u^*(x)$, then $\min\{\mathbb I_\alpha u^*(x),g_2(x)\}>u^*(x)\geq g_1(x).$ It would follow from (3.6) that $u^*(x)>u^*(x)$, which is impossible. Therefore, $\mathbb I_\alpha u^*(x)\leq u^*(x).$

		\item[(c)] Suppose that $g_1(x)<u^*(x)\leq g_2(x)$. It follows that $g_2(x)>g_1(x)$, together with (3.6) gives $u^*(x)=\min\{\mathbb I_\alpha u^*(x),g_2(x)\}.$
		Therefore, $u^*(x)\leq\mathbb I_\alpha u^*(x).$
		\end{itemize}
	\end{proof}
	
	For $u^*$ satisfying (3.6), define
	\begin{equation}
		D^1_{u^*}:=\{x\in D^c:g_1(x)\geq u^*(x)\}; \ \ \ \ \ D^2_{u^*}:=\{x\in D^c:g_2(x)\leq u^*(x)\}. \tag{4.2}
	\end{equation}
	Define the corresponding stopping times by
	\begin{equation}
		\theta^1_{u^*}:=\inf\{t|X_t\in D^1_{u^*}\}; \ \ \ \ \ \theta^2_{u^*}:=\inf\{t|X_t\in D^2_{u^*}\}. \tag{4.3}
	\end{equation}
	Here $\inf\varnothing:=+\infty$. We next verify that the two expressions in (4.3) are well defined as stopping times. Since $X$ is a countable discrete space, $D_{u^*}^i\in\mathcal{B}(X)$ for $i=1,2$. Moreover, for every $n\in\mathbb{N}_0$,
	\begin{equation*}
		\{\theta_{u^*}^i\leq n\}
		=\bigcup_{t=0}^{n}\{X_t\in D_{u^*}^i\}\in\mathcal{G}_n,
		\qquad i=1,2,
	\end{equation*}
	because $X_t$ is $\mathcal{G}_t$-measurable and $\mathcal{G}_t\subseteq\mathcal{G}_n$ for $t\leq n$. Hence $\theta_{u^*}^1,\theta_{u^*}^2\in\Theta$.

	By the measurable selection theorem in \autocite{Nowak_1985}, selecting controls on $D^c$ and extending them arbitrarily to $D$, there exists a pair $(\pi^1_{u^*},\pi^2_{u^*})\in\Pi_s^1\times\Pi_s^2$ such that
	\begin{align}
		\mathbb{I}_\alpha u^*(x) &= r(x,\pi^1_{u^*},\pi^2_{u^*}) + \alpha(x)\sum_{y\in X} u^*(y)p(y|x,\pi^1_{u^*},\pi^2_{u^*}) \notag \\
		&=\sup_{\lambda\in \mathcal{P}(A(x))}\bigg[r(x,\lambda,\pi^2_{u^*}) + \alpha(x)\sum_{y\in X} u^*(y)p(y|x,\lambda,\pi^2_{u^*})\bigg]  \notag \\ 
		&=\inf_{\mu\in \mathcal{P}(B(x))}\bigg[r(x,\pi^1_{u^*},\mu) + \alpha(x)\sum_{y\in X} u^*(y)p(y|x,\pi^1_{u^*},\mu)\bigg]  \ \ \ \forall x\in D^c.\tag{4.4}
	\end{align}
	The next theorem is the principal existence result of this section.
	\begin{theorem}
		Under Assumptions 3.1 and 3.3, let $u^*$ be the solution of (3.6) given by Theorem 3.4, and let $(\pi^1_{u^*},\theta^1_{u^*};\pi^2_{u^*},\theta^2_{u^*})$ be defined by (4.2)--(4.4). Then:
		\begin{itemize}
		\item[(a)] The value of the game exists and $J(x)=u^*(x)$ for all $x\in D^c$.
		\item[(b)] The strategy pair $(\pi^1_{u^*},\theta^1_{u^*};\pi^2_{u^*},\theta^2_{u^*})$ is a Nash equilibrium.
		\end{itemize}
	\end{theorem}

	\begin{proof}
		\begin{itemize}
		\item[(a)] Fix $x\in D^c$. Let $(\pi^2,\theta^2)\in\Pi^2\times\Theta$ be arbitrary. By the last equality in (4.4) and by (2.4)--(2.5), for every $s\in\mathbb N_0$, on $\{s<\sigma_D\}$,
		\begin{align}
			\mathbb I_\alpha u^*(X_s)
			&\leq\int_{B(X_s)}
			\bigg[r(X_s,\pi_{u^*}^1,b)
			+\alpha(X_s)\sum_{y\in X}u^*(y)p(y\mid X_s,\pi_{u^*}^1,b)\bigg]
			\pi_s^2(db\mid Z_s)\notag\\
			&=E_x^{\pi_{u^*}^1,\pi^2}\!\left[
			r(X_s,A_s,B_s)+\alpha(X_s)u^*(X_{s+1})
			\mid\mathcal G_s\right],
			\quad P_x^{\pi_{u^*}^1,\pi^2}\text{-a.s.} \tag{4.5}
		\end{align}

		We next use the definition of $\theta_{u^*}^1$ to determine which conclusion of Lemma 4.1 applies before player 1 stops or the system fails. On $\{s<\theta_{u^*}^1\wedge\sigma_D\}$,
		\[
			g_1(X_s)<u^*(X_s)\leq g_2(X_s).
		\]
		Thus, on $\{s<\theta_{u^*}^1\wedge\sigma_D\}$, Lemma 4.1(c) and (4.5) give
		\begin{align}
			u^*(X_s)
			\leq E_x^{\pi_{u^*}^1,\pi^2}\!\left[
			r(X_s,A_s,B_s)+\alpha(X_s)u^*(X_{s+1})
			\mid\mathcal G_s\right]
			\quad P_x^{\pi_{u^*}^1,\pi^2}\text{-a.s.} \tag{4.6}
		\end{align}
		
		To relate (4.6) to the discounted criterion, fix a deterministic $T\in\mathbb N$. We first sum the discounted inequalities over $s=0,\ldots,T-1$, retaining only the terms before either player stops or the system fails. The passage to $T\to\infty$ is justified below. For $0\leq s\leq T-1$, the product $\left(\prod_{k=0}^{s-1}\alpha(X_k)\right)\allowbreak I_{\{s<T\wedge\theta_{u^*}^1\wedge\theta^2\wedge\sigma_D\}}$
		is $\mathcal G_s$-measurable. Hence,
		\begin{align*}
			&E_x^{\pi_{u^*}^1,\pi^2}\!\bigg[
			\left(\prod_{k=0}^{s-1}\alpha(X_k)\right)I_{\{s<T\wedge\theta_{u^*}^1\wedge\theta^2\wedge\sigma_D\}}
			E_x^{\pi_{u^*}^1,\pi^2}\!\left[
			r(X_s,A_s,B_s)+\alpha(X_s)u^*(X_{s+1})
			\mid\mathcal G_s\right]\bigg]\\
			&\qquad=
			E_x^{\pi_{u^*}^1,\pi^2}\!\bigg[
			\left(\prod_{k=0}^{s-1}\alpha(X_k)\right)I_{\{s<T\wedge\theta_{u^*}^1\wedge\theta^2\wedge\sigma_D\}}
			\left(r(X_s,A_s,B_s)+\alpha(X_s)u^*(X_{s+1})\right)\bigg].
		\end{align*}
		Multiplying (4.6) by $\left(\prod_{k=0}^{s-1}\alpha(X_k)\right)I_{\{s<T\wedge\theta_{u^*}^1\wedge\theta^2\wedge\sigma_D\}},$
		taking expectations, and summing over $s=0,\ldots,T-1$ gives
		\begin{align}
			E_x^{\pi_{u^*}^1,\pi^2}\!\left[
			\sum_{s=0}^{T\wedge\theta_{u^*}^1\wedge\theta^2\wedge\sigma_D-1} \left(\prod_{k=0}^{s-1} \alpha(X_k)\right)u^*(X_s)\right]
			&\leq E_x^{\pi_{u^*}^1,\pi^2}\!\left[
			\sum_{s=0}^{T\wedge\theta_{u^*}^1\wedge\theta^2\wedge\sigma_D-1}
			\left(\prod_{k=0}^{s-1}\alpha(X_k)\right)r(X_s,A_s,B_s)\right]\notag\\
			&+ E_x^{\pi_{u^*}^1,\pi^2}\!\left[
			\sum_{s=0}^{T\wedge\theta_{u^*}^1\wedge\theta^2\wedge\sigma_D-1}
			\left(\prod_{k=0}^{s}\alpha(X_k)\right)u^*(X_{s+1})\right]. \tag{4.7}
		\end{align}
		Subtracting the two sums involving $u^*$ in (4.7) gives
		\begin{align*}
			&\sum_{s=0}^{T\wedge\theta_{u^*}^1\wedge\theta^2\wedge\sigma_D-1}
			\left(\prod_{k=0}^{s-1}\alpha(X_k)\right)u^*(X_s)
			-\sum_{s=0}^{T\wedge\theta_{u^*}^1\wedge\theta^2\wedge\sigma_D-1}
			\left(\prod_{k=0}^{s}\alpha(X_k)\right)u^*(X_{s+1})\\
			&\qquad=u^*(x)
			-\left(\prod_{k=0}^{T\wedge\theta_{u^*}^1\wedge\theta^2\wedge\sigma_D-1}\alpha(X_k)\right)
			u^*(X_{T\wedge\theta_{u^*}^1\wedge\theta^2\wedge\sigma_D}).
		\end{align*}
		Substituting this identity into (4.7) gives
		\begin{align}
			u^*(x)
			&\leq E_x^{\pi_{u^*}^1,\pi^2}\!\left[
			\sum_{s=0}^{T\wedge\theta_{u^*}^1\wedge\theta^2\wedge\sigma_D-1}
			\left(\prod_{k=0}^{s-1}\alpha(X_k)\right)r(X_s,A_s,B_s)\right]\notag\\
			&\quad+
			E_x^{\pi_{u^*}^1,\pi^2}\!\left[
			\left(\prod_{k=0}^{T\wedge\theta_{u^*}^1\wedge\theta^2\wedge\sigma_D-1}\alpha(X_k)\right)
			u^*(X_{T\wedge\theta_{u^*}^1\wedge\theta^2\wedge\sigma_D})\right]. \tag{4.8}\label{eq:verify-truncated-1}
		\end{align}

		Before system failure, the definition of $D_{u^*}^1$ and $u^*\geq g_1$ on $D^c$ give
		\begin{equation}
			u^*(X_{\theta_{u^*}^1})=g_1(X_{\theta_{u^*}^1})
			\quad\text{on }\{\theta_{u^*}^1<\sigma_D\}. \tag{4.9}\label{eq:verify-stop-1}
		\end{equation}
		On $\{\theta^2<\theta_{u^*}^1,\ \theta^2<\sigma_D\}$, one has
		\begin{equation}
			u^*(X_{\theta^2})\leq g_2(X_{\theta^2}). \tag{4.10}\label{eq:verify-opponent-stop-2}
		\end{equation}

		If system failure occurs first or simultaneously with a player's stopping time, then
		\begin{equation}
			u^*(X_{\sigma_D})=0
			\quad\text{on }\{\sigma_D\leq\theta_{u^*}^1\wedge\theta^2,\ \sigma_D<\infty\}. \tag{4.11}\label{eq:verify-failure-1}
		\end{equation}
		
		Letting $T\to\infty$ in (4.8), we obtain by the dominated convergence theorem and \eqref{eq:verify-stop-1}--\eqref{eq:verify-failure-1}
		\begin{align}
			u^*(x)
			\leq E_x^{\pi_{u^*}^1,\pi^2}\bigg[&
			\sum_{s=0}^{\theta_{u^*}^1\wedge\theta^2\wedge\sigma_D-1}
			\left(\prod_{k=0}^{s-1}\alpha(X_k)\right)r(X_s,A_s,B_s)
			+\left(\prod_{k=0}^{\theta_{u^*}^1-1}\alpha(X_k)\right)u^*(X_{\theta_{u^*}^1})
			I_{\{\theta_{u^*}^1\leq\theta^2,\,\theta_{u^*}^1<\sigma_D\}}\notag\\
			&\mathrel{\phantom{\leq}}{+}\left(\prod_{k=0}^{\theta^2-1}\alpha(X_k)\right)u^*(X_{\theta^2})
			I_{\{\theta^2<\theta_{u^*}^1,\,\theta^2<\sigma_D\}}
			\bigg]\quad\text{by (4.8)}\notag\displaybreak[1]\\
			=E_x^{\pi_{u^*}^1,\pi^2}\bigg[&
			\sum_{s=0}^{\theta_{u^*}^1\wedge\theta^2\wedge\sigma_D-1}
			\left(\prod_{k=0}^{s-1}\alpha(X_k)\right)r(X_s,A_s,B_s)
			+\left(\prod_{k=0}^{\theta_{u^*}^1-1}\alpha(X_k)\right)g_1(X_{\theta_{u^*}^1})
			I_{\{\theta_{u^*}^1\leq\theta^2,\,\theta_{u^*}^1<\sigma_D\}}\notag\\
			&\mathrel{\phantom{=}}{+}\left(\prod_{k=0}^{\theta^2-1}\alpha(X_k)\right)u^*(X_{\theta^2})
			I_{\{\theta^2<\theta_{u^*}^1,\,\theta^2<\sigma_D\}}
			\bigg]\quad\text{by (4.9)}\notag\displaybreak[1]\\
			\leq E_x^{\pi_{u^*}^1,\pi^2}\bigg[&
			\sum_{s=0}^{\theta_{u^*}^1\wedge\theta^2\wedge\sigma_D-1}
			\left(\prod_{k=0}^{s-1}\alpha(X_k)\right)r(X_s,A_s,B_s)
			+\left(\prod_{k=0}^{\theta_{u^*}^1-1}\alpha(X_k)\right)g_1(X_{\theta_{u^*}^1})
			I_{\{\theta_{u^*}^1\leq\theta^2,\,\theta_{u^*}^1<\sigma_D\}}\notag\\
			&\mathrel{\phantom{\leq}}{+}\left(\prod_{k=0}^{\theta^2-1}\alpha(X_k)\right)g_2(X_{\theta^2})
			I_{\{\theta^2<\theta_{u^*}^1,\,\theta^2<\sigma_D\}}
			\bigg]\quad\text{by (4.10)}\notag\\
			=J(x,\pi_{u^*}^1,&\theta_{u^*}^1,\pi^2,\theta^2)\qquad\text{by (2.6)}. \tag{4.12}\label{eq:verify-criterion-1}
		\end{align}
		For arbitrary $(\pi^1,\theta^1)\in\Pi^1\times\Theta$, on  $\{s<\theta_{u^*}^2\wedge\sigma_D\}$, using (4.4) and Lemma 4.1(b) gives
		\begin{align}
			u^*(X_s)
			\geq E_x^{\pi^1,\pi_{u^*}^2}\!\left[
			r(X_s,A_s,B_s)+\alpha(X_s)u^*(X_{s+1})
			\mid\mathcal G_s\right]
			\quad P_x^{\pi^1,\pi_{u^*}^2}\text{-a.s.} \tag{4.13}
		\end{align}
		Proceeding as in (4.6)--(4.8), for every $T\in\mathbb N$,
		\begin{align}
			u^*(x)
			&\geq E_x^{\pi^1,\pi_{u^*}^2}\!\left[
			\sum_{s=0}^{T\wedge\theta^1\wedge\theta_{u^*}^2\wedge\sigma_D-1}
			\left(\prod_{k=0}^{s-1}\alpha(X_k)\right)r(X_s,A_s,B_s)\right]\notag\\
			&\quad+
			E_x^{\pi^1,\pi_{u^*}^2}\!\left[
			\left(\prod_{k=0}^{T\wedge\theta^1\wedge\theta_{u^*}^2\wedge\sigma_D-1}\alpha(X_k)\right)
			u^*(X_{T\wedge\theta^1\wedge\theta_{u^*}^2\wedge\sigma_D})\right]. \tag{4.14}\label{eq:verify-truncated-2}
		\end{align}
		Since $u^*\geq g_1$ on $D^c$,
		\begin{equation}
			u^*(X_{\theta^1})\geq g_1(X_{\theta^1})
			\quad\text{on }\{\theta^1<\sigma_D\}. \tag{4.15}\label{eq:verify-opponent-stop-1}
		\end{equation}
		Moreover,
		\begin{equation}
			u^*(X_{\theta_{u^*}^2})\geq g_2(X_{\theta_{u^*}^2})
			\quad\text{on }\{\theta_{u^*}^2<\sigma_D\}. \tag{4.16}\label{eq:verify-stop-2}
		\end{equation}
		Letting $T\to\infty$ in (4.14) and applying the same argument as above, with (4.15)--(4.16) and the zero extension on $D$, gives
		\begin{equation}
			u^*(x)\geq J(x,\pi^1,\theta^1,\pi^2_{u^*},\theta^2_{u^*}). \tag{4.17}\label{eq:verify-criterion-2}
		\end{equation}

		We now combine (4.12) and (4.17) to prove the existence of the value. Since (4.12) holds for every strategy of player 2 and (4.17) holds for every strategy of player 1,
		\begin{align}
			u^*(x)&\leq\inf_{(\pi^2,\theta^2)\in\Pi^2\times\Theta}J(x,\pi_{u^*}^1,\theta_{u^*}^1,\pi^2,\theta^2)
			\leq J^-(x)\leq J^+(x) \notag\\
			&\leq\sup_{(\pi^1,\theta^1)\in\Pi^1\times\Theta}J(x,\pi^1,\theta^1,\pi_{u^*}^2,\theta_{u^*}^2)\leq u^*(x). \tag{4.18}
		\end{align}
		Thus all inequalities in (4.18) are equalities. In particular,
		\[
			J^-(x)=J^+(x)=u^*(x),
		\]
		so the value exists and $J(x)=u^*(x)$, proving part (a).

		\item[(b)] By (4.12), (4.17), and part (a), for every $(\pi^1,\theta^1)\in\Pi^1\times\Theta$ and $(\pi^2,\theta^2)\in\Pi^2\times\Theta$,
		\begin{align*}
			J(x,\pi^1,\theta^1,\pi_{u^*}^2,\theta_{u^*}^2)
			\leq J(x,\pi_{u^*}^1,\theta_{u^*}^1,\pi_{u^*}^2,\theta_{u^*}^2)
			=u^*(x)
			\leq J(x,\pi_{u^*}^1,\theta_{u^*}^1,\pi^2,\theta^2).
		\end{align*}
		Therefore, $(\pi_{u^*}^1,\theta_{u^*}^1;\pi_{u^*}^2,\theta_{u^*}^2)$ is a Nash equilibrium. This proves part (b).
		\end{itemize}
	\end{proof}
	
	\begin{remark}
		Theorem 4.2 combines assumptions that are treated separately in the related control-and-stopping literature. Compared with the \autocite{Ghosh_2007}, which assumes bounded rewards and strictly ordered terminal rewards in its main results and discusses an extension through weighted norms, this paper allows the discount factor to depend on the current state and gives explicit conditions for possibly unbounded running and terminal rewards. Compared with \autocite{Chen_Guo_Guo_Jia_2025}, the present model uses state-dependent discount factors rather than a constant discount rate and does not require bounded rewards. 
	\end{remark}

	\section{Computation of the value and equilibrium}\label{s:5}
	This section develops algorithms for approximating the value and constructing an $\varepsilon$-Nash equilibrium. We prove approximation of the value $J$ using the sequence $\{\hat u_n\}$ defined in (3.5). Since $\hat u_n\leq J$, the approximation error at state $x$ after $n$ iterations is $J(x)-\hat u_n(x)$. We bound this error using the difference between two successive iterates. Finite action sets are assumed for the linear-programming computations, not for the error bounds.
	\subsection{Computation of the value}
	Theorems 3.4 and 4.2 show that $\{\hat u_n\}$ in (3.5) increases to the value $J$. Following the contraction argument and value-iteration method in Puterman \autocite[Sections 6.2.3 and 6.3.2]{Puterman_1994}, we use the contraction inequality (3.8) to bound the error from the difference $\hat u_n-\hat u_{n-1}$. Set $\kappa:=\bar\alpha m<1$, as in the proof of Theorem 3.4.
	\begin{theorem}
		Under Assumptions 3.1 and 3.3, the sequence $\{\hat u_n\}$ in (3.5) has the following properties.
		\begin{itemize}
		\item[(a)] $J=u^*$, where $u^*$ is the fixed point of $\mathbb T_{\alpha}$ in (3.4). For $x\in D^c$ and $n\geq0$,
		\[
			g_1(x)\leq\hat u_n(x)\leq\hat u_{n+1}(x)
			\leq J(x)\leq\max\{g_1(x),g_2(x)\}.
		\]
		Moreover, $\{\hat u_n\}$ converges to $J$ in $\lVert \cdot\rVert_{V_0}$.
		\item[(b)] For every $n\geq1$,
		\begin{equation}
			\lVert J-\hat u_n\rVert_{V_0}
			\leq\frac{\kappa}{1-\kappa}\lVert \hat u_n-\hat u_{n-1}\rVert_{V_0}
			\leq\frac{\kappa^n}{1-\kappa}\lVert \hat u_1-\hat u_0\rVert_{V_0}. \tag{5.1}
		\end{equation}
		\end{itemize}
	\end{theorem}
	\begin{proof}
		\begin{itemize}
		\item[(a)] By Theorem 4.2, $J=u^*$. Theorem 3.4 gives $\hat u_n(x)\leq\hat u_{n+1}(x)$ and weighted-norm convergence of $\{\hat u_n\}$. Since $\hat u_0\leq J$ and $\mathbb T_\alpha J=J$, the monotonicity of $\mathbb T_\alpha$ gives
		\[
			\hat u_n \leq \hat u_{n+1}=\mathbb T_\alpha^{n+1}\hat u_0
			\leq\mathbb T_\alpha^{n+1}J=J.
		\]
		By (3.5)--(3.6), we have $g_1(x)\leq\hat u_n(x)\leq\hat u_{n+1}(x) \leq J(x)\leq\max\{g_1(x),g_2(x)\}$ on $D^c$. This proves (a).

		\item[(b)] The contraction property (3.8) and $\mathbb T_\alpha J=J$ give (5.1). By part (a) and the definition of the weighted norm,
		\begin{equation}
			0\leq J(x)-\hat u_n(x)
			\leq\frac{\kappa}{1-\kappa}\lVert \hat u_n-\hat u_{n-1}\rVert_{V_0}V_0(x),
			\qquad x\in D^c,\quad n\geq1. \tag{5.2}
		\end{equation}
		\end{itemize}
	\end{proof}

	To compute $\hat u_n(x)$ for $x\in D^c$, it suffices to evaluate $\mathbb I_\alpha u(x)$ with $u=\hat u_{n-1}$. When the action sets are finite, this quantity is the value of a zero-sum matrix game and can be computed by linear programming.

	The primal program is
	\begin{equation} \tag{5.3}
		\begin{array}{rl}
			& \max v \\[3pt]
			& \left\{
			\begin{aligned}
				& v \leq \sum_{a \in A(x)} \mathbb{I}_\alpha u(x,a,b)\phi(a), \quad b \in B(x) \\
				& \sum_{a \in A(x)} \phi(a) = 1 \\
				& \phi(a) \geq 0, \, a \in A(x),
			\end{aligned}
			\right.
		\end{array}
	\end{equation}
	where $\mathbb{I}_\alpha u(x,a,b)$ is defined as $\mathbb{I}_\alpha u(x,a,b):= r(x,a,b) + \alpha(x)\sum_{y\in X} u(y)p(y|x,a,b)$.
	
	The dual program is
	\begin{equation} \tag{5.4}
		\begin{array}{rl}
			& \min z \\[3pt]
			& \left\{
			\begin{aligned}
				& z \geq \sum_{b \in B(x)} \mathbb{I}_\alpha u(x,a,b)\psi(b), \quad a \in A(x) \\
				& \sum_{b \in B(x)} \psi(b) = 1 \\
				& \psi(b) \geq 0, \, b \in B(x).
			\end{aligned}
			\right.
		\end{array}
	\end{equation}
	When $A(x)$ and $B(x)$ are finite, both programs admit solutions. Denote optimal solutions of (5.3) and (5.4) by $\phi_u(\cdot\mid x)$ and $\psi_u(\cdot\mid x)$, respectively. They form a saddle point for the one-stage zero-sum game with payoff $\mathbb{I}_\alpha u(x,a,b)$; see, for example, \autocite[Section 2.4]{Barron_2013}. Consequently,
	\begin{equation}
		\mathbb{I}_\alpha u(x)=\sum_{a \in A(x)} \sum_{b \in B(x)}\mathbb{I}_\alpha u(x,a,b)\phi_u(a|x)\psi_u(b|x). \tag{5.5}
	\end{equation}
	For $x\in D^c$, applying (5.5) with $u=\hat u_{n-1}$ gives $\mathbb I_\alpha\hat u_{n-1}(x)$ and hence
	\begin{equation}
		\hat u_n(x)=\max\{\min\{\mathbb I_\alpha\hat u_{n-1}(x),g_2(x)\},g_1(x)\},
		\qquad n\geq1. \tag{5.6}
	\end{equation}
	
	The bound in (5.2) contains $V_0(x)$, when $V_0$ is unbounded,
	it does not guarantee that
	\[
	0\leq J(x)-\hat u_n(x)\leq\delta,\qquad x\in D^c,
	\]
	for a prescribed $\delta>0$ and a single iteration index $n$.
	After explaining how to compute $\hat u_n$ by linear programming,
	we next give an additional condition under which this inequality holds after finite iterations. We then use it to construct an $\varepsilon$-Nash equilibrium. For a bounded function $f$ on $X$, we write $\lVert f\rVert_\infty:=\sup_{x\in X}|f(x)|$.
	\begin{theorem}
		Suppose that Assumptions 3.1 and 3.3 hold and that
		\[
			X_*:=\{x\in D^c:g_1(x)\neq g_2(x)\}
		\]
		is finite. Then the following statements hold.
		\begin{itemize}
		\item[(a)] $\hat u_n(x)=J(x)$ for all $x\in D^c\setminus X_*$ and $n\geq0$.
		\item[(b)] For every $n\geq1$ and $x\in D^c$,
		\[
			0\leq J(x)-\hat u_n(x)
			\leq\frac{\bar\alpha}{1-\bar\alpha}
			\max_{y\in X_*}|\hat u_n(y)-\hat u_{n-1}(y)|.
		\]
		\item[(c)] For every $\delta>0$, there exists an integer $n\geq1$
		such that
		\[
		\qquad
		\frac{\bar\alpha}{1-\bar\alpha}
		\max_{x\in X_*}
		|\hat u_n(x)-\hat u_{n-1}(x)|
		\leq\delta.
		\]
		When $X_*=\varnothing$, the maxima over $X_*$ below
		are defined to be zero, and we have $\hat u_0=J$.
		\end{itemize}
	\end{theorem}
	\begin{proof}
		\begin{itemize}
		\item[(a)] If $x\in D^c\setminus X_*$, then $g_1(x)=g_2(x)$. From (3.4), (3.5) and (3.6) we have
		\[
			\hat u_n(x)=g_1(x)=g_2(x)=J(x),\qquad n\geq0.
		\]
		This proves (a).

		\item[(b)] For $u,v\in M_{V_0}(X)$, if $u-v$ is bounded, by (3.3), for every $x\in D^c$,
		\begin{align*}
			|\mathbb I_\alpha u(x)-\mathbb I_\alpha v(x)|
			&\leq
			\sup_{\lambda\in\mathcal P(A(x))}
			\sup_{\mu\in\mathcal P(B(x))}
			\bigg|
			\int_{A(x)}\int_{B(x)}
			\bigl(\mathbb I_\alpha u(x,a,b)
			-\mathbb I_\alpha v(x,a,b)\bigr)
			\mu(db)\lambda(da)
			\bigg|\\
			&\leq
			\alpha(x)
			\sup_{\lambda\in\mathcal P(A(x))}
			\sup_{\mu\in\mathcal P(B(x))}
			\int_{A(x)}\int_{B(x)}
			\sum_{y\in X}|u(y)-v(y)|p(y\mid x,a,b)
			\mu(db)\lambda(da)\\
			&\leq
			\alpha(x)\sup_{a\in A(x),\,b\in B(x)}
			\sum_{y\in X}|u(y)-v(y)|p(y\mid x,a,b)\\
			&\leq \alpha(x)\lVert u-v\rVert_\infty
			\leq\bar\alpha\lVert u-v\rVert_\infty.
		\end{align*}
		Using (3.7) yields
		\begin{align}
			\lVert \mathbb T_\alpha u-\mathbb T_\alpha v\rVert_\infty
			\leq\sup_{x\in D^c}
			|\mathbb I_\alpha u(x)-\mathbb I_\alpha v(x)|
			\leq\bar\alpha\lVert u-v\rVert_\infty. \tag{5.7}\label{eq:uniform-contraction}
		\end{align}
		By part (a), $\lVert J-\hat u_{n-1}\rVert_\infty = \max_{x\in X_*}|J(x)-\hat u_{n-1}(x)|<\infty, \forall n\geq1.$ Therefore, by \eqref{eq:uniform-contraction},
		\begin{align*}
			\lVert J-\hat u_n\rVert_\infty
			=\lVert \mathbb T_\alpha J-\mathbb T_\alpha \hat u_{n-1}\rVert_\infty
			\leq\bar\alpha\lVert J-\hat u_{n-1}\rVert_\infty
			\leq\bar\alpha\lVert J-\hat u_n\rVert_\infty
			+\bar\alpha\lVert \hat u_n-\hat u_{n-1}\rVert_\infty.
		\end{align*}
		Rearranging and using part (a) and Theorem 5.1(a) gives
		\begin{align*}
			0\leq J(x)-\hat u_n(x)
			\leq\lVert J-\hat u_n\rVert_\infty
			\leq\frac{\bar\alpha}{1-\bar\alpha}
			\lVert \hat u_n-\hat u_{n-1}\rVert_\infty
			=\frac{\bar\alpha}{1-\bar\alpha}
			\max_{y\in X_*}|\hat u_n(y)-\hat u_{n-1}(y)|.
		\end{align*}
		This proves part (b).
		
		\item[(c)] Since $X_*$ is finite, from Theorem~5.1(a),
         $\lim_{n\to\infty}
         \max_{x\in X_*}|\hat u_n(x)-\hat u_{n-1}(x)|=0.$
		Consequently, for every $\delta>0$, there exists an integer
		$n_0\geq1$ such that
		$\frac{\bar\alpha}{1-\bar\alpha}
		\max_{x\in X_*}|\hat u_n(x)-\hat u_{n-1}(x)|
		\leq\delta, \ n\geq n_0.$ If $X_*=\varnothing$, part (a) gives $\hat u_0=J$ directly.
		\end{itemize}
	\end{proof}

	Theorem 5.2(c) gives the following stopping rule, while Theorem 5.2(b) gives the corresponding error bound. A finite $X_*$ ensures that the rule is attained.

	\noindent\begin{minipage}{\linewidth}
	\textbf{Algorithm \RN{1} (Approximation of the value of the game):}\\
	Step 1 (Initialization): Input $\delta>0$, set $n=0$, and take $\hat u_0$ from (3.5). Form the finite set $X_*$. If $X_*=\varnothing$, return $J_\delta=\hat u_0$.\\
	Step 2 (Value update): Replace $n$ by $n+1$. Compute $\hat u_n(x)$ for $x\in X_*$ using (5.3)--(5.6). Set $\hat u_n(x)=g_1(x)$ on $D^c\setminus X_*$.\\
	Step 3 (Stopping rule): If
	\[
		\max_{x\in X_*}|\hat u_n(x)-\hat u_{n-1}(x)|
		\leq\frac{1-\bar\alpha}{\bar\alpha}\delta,
	\]
	go to Step 4. Otherwise, return to Step 2.\\
	Step 4 (Value approximation): Return $J_\delta=\hat u_n$. Then
	$0\leq J(x)-J_\delta(x)\leq\delta$ for every $x\in D^c$.
	\end{minipage}
	\subsection{Computation of an \texorpdfstring{$\varepsilon$}{epsilon}-Nash equilibrium}
	The $\delta$ above denotes the value-approximation accuracy in Algorithm \RN{1}, whereas $\varepsilon$ below denotes the error in the Nash equilibrium inequalities. For a given $\varepsilon>0$, set
	\[
		\delta_\varepsilon:=\frac{\bar\alpha(1-\bar\alpha)}{6}\varepsilon.
	\]
	
	Run Algorithm \RN{1} with accuracy $\delta_\varepsilon$. Its stopping condition becomes
	\[
		\max_{x\in X_*}|\hat u_n(x)-\hat u_{n-1}(x)|
		\leq\frac{(1-\bar\alpha)^2}{6}\varepsilon.
	\]
	If $X_*=\varnothing$, take $\tilde n=0$; otherwise, let $\tilde n$ be the stopping index. Theorem 5.2 then gives
	\begin{align}
		0\leq J(x)-J_{\delta_\varepsilon}(x)
		\leq\delta_\varepsilon
		=\frac{\bar\alpha(1-\bar\alpha)}{6}\varepsilon,
		\quad x\in D^c,
		\qquad J_{\delta_\varepsilon}(x):=\hat u_{\tilde n}(x). \tag{5.8}
	\end{align}
	Then, by Theorems 5.1 and 5.2,
	\begin{equation}
		\begin{cases}
			g_1(x)\leq J_{\delta_\varepsilon}(x)\leq\max\{g_1(x),g_2(x)\},&x\in D^c,\\
			J_{\delta_\varepsilon}(x)=g_1(x)=g_2(x),&x\in D^c\setminus X_*.
		\end{cases}\tag{5.9}
	\end{equation}
	Using $J_{\delta_\varepsilon}$, first define
	\begin{align}
		D_\varepsilon^1:=\{x\in D^c:g_1(x)-J_{\delta_\varepsilon}(x)\geq -\frac{\varepsilon}{4}\}; \ \ D_\varepsilon^2:=\{x\in D^c:g_2(x)-J_{\delta_\varepsilon}(x)\leq \frac{\varepsilon}{4}\}. \tag{5.10}
	\end{align}
	The associated stopping times are
	\begin{align}
		\theta_\varepsilon^1:=\inf\{t\geq0:X_t\in D_\varepsilon^1\}; \ \ \ \theta_\varepsilon^2:=\inf\{t\geq0:X_t\in D_\varepsilon^2\}, \tag{5.11}
	\end{align}
	where $\inf\varnothing:=+\infty$. The stopping times in (5.11) are well defined. Indeed, for $k=1,2$, we have $\theta_\varepsilon^k\in\Theta$, because
	\[
		\{\theta_\varepsilon^k\leq n\}
		=\bigcup_{t=0}^{n}\{X_t\in D_\varepsilon^k\}
		\in\mathcal G_n, \ \ \forall n\in\mathbb N_0.
	\]
	Solving (5.3) and (5.4) for $x\in D^c$ with $u=J_{\delta_\varepsilon}$, and extending the controls arbitrarily to $D$, as in \autocite[Section 2.4]{Barron_2013}, yields policies $\pi_\varepsilon^1(\cdot\mid x)$ and $\pi_\varepsilon^2(\cdot\mid x)$ such that $(\pi_\varepsilon^1,\pi_\varepsilon^2)\in\Pi_s^1\times\Pi_s^2$ and
	\begin{align}
		\mathbb{I}_\alpha J_{\delta_\varepsilon}(x) &= r(x,\pi^1_\varepsilon,\pi^2_\varepsilon) + \alpha(x)\sum_{y\in X} J_{\delta_\varepsilon}(y)p(y|x,\pi^1_\varepsilon,\pi^2_\varepsilon) \notag \\
		&=\sup_{\lambda\in \mathcal{P}(A(x))}\bigg[r(x,\lambda,\pi^2_\varepsilon) + \alpha(x)\sum_{y\in X} J_{\delta_\varepsilon}(y)p(y|x,\lambda,\pi^2_\varepsilon)\bigg]  \notag \\ 
		&=\inf_{\mu\in \mathcal{P}(B(x))}\bigg[r(x,\pi^1_\varepsilon,\mu) + \alpha(x)\sum_{y\in X} J_{\delta_\varepsilon}(y)p(y|x,\pi^1_\varepsilon,\mu)\bigg],  \ \ \ \forall x\in D^c.\tag{5.12}
	\end{align}
	
	The next theorem justifies this construction.
	\begin{theorem}
		Under Assumptions 3.1 and 3.3, suppose that $X_*:=\{x\in D^c:g_1(x)\neq g_2(x)\}$ is finite and that the action sets $A(x)$ and $B(x)$ are finite for every $x\in D^c$. Then, for every $\varepsilon>0$, $(\pi_\varepsilon^1,\theta_\varepsilon^1;\pi_\varepsilon^2,\theta_\varepsilon^2)$ is an $\varepsilon$-Nash equilibrium.
	\end{theorem}
	\begin{proof}
		By Theorem~5.1, $J=u^*$ on $D^c$, so
		\begin{align}
			J(x)=\max\{\min\{\mathbb I_\alpha J(x),g_2(x)\},g_1(x)\},
			\qquad x\in D^c. \tag{5.13}
		\end{align}
		Put $\Delta_\varepsilon(x):=J(x)-J_{\delta_\varepsilon}(x).$
		By (5.8), $0\leq\Delta_\varepsilon(x)\leq\delta_\varepsilon$ for every $x\in D^c$. Define $d_\varepsilon(x):=\mathbb I_\alpha J(x)-\mathbb I_\alpha J_{\delta_\varepsilon}(x).$
		Since $0\leq\Delta_\varepsilon\leq\delta_\varepsilon$, for every $(x,a,b)\in K$, $0\leq\alpha(x)\sum_{y\in X}\Delta_\varepsilon(y)p(y\mid x,a,b) \leq\alpha(x)\delta_\varepsilon.$
		The monotonicity of $\mathbb I_\alpha$ and (4.1) therefore give
		\begin{align}
			0\leq d_\varepsilon(x)\leq\alpha(x)\delta_\varepsilon\leq\bar\alpha\delta_\varepsilon,
			\qquad x\in D^c. \tag{5.14}
		\end{align}

		On the one hand, we show that $J(x,\pi_\varepsilon^1,\theta_\varepsilon^1,\pi^2,\theta^2)$ is at least $J(x)-\varepsilon$.
		Fix $x\in D^c$ and $(\pi^2,\theta^2)\in\Pi^2\times\Theta$ and work under $P_x^{\pi_\varepsilon^1,\pi^2}$.
		For $s<\theta_\varepsilon^1\wedge\sigma_D$, the definition of $D_\varepsilon^1$ and (5.8) yield
		\[
			g_1(X_s)<J_{\delta_\varepsilon}(X_s)-\frac{\varepsilon}{4}<J(X_s).
		\]
		Since $J(X_s)>g_1(X_s)$, (3.6) also gives $J(X_s)\leq g_2(X_s)$. Lemma~4.1(c) therefore gives $J(X_s)\leq\mathbb I_\alpha J(X_s)$. Combining this inequality with (5.12), (5.14), and conditional expectation with respect to $\mathcal G_s=\sigma( Z_s)$ gives
		\begin{align}
			J_{\delta_\varepsilon}(X_s)
			&\leq E_x^{\pi_\varepsilon^1,\pi^2}\!\left[
			r(X_s,A_s,B_s)+\alpha(X_s)J_{\delta_\varepsilon}(X_{s+1})
			\mid\mathcal G_s\right]
			+d_\varepsilon(X_s)-\Delta_\varepsilon(X_s)\notag\\
			&\leq E_x^{\pi_\varepsilon^1,\pi^2}\!\left[
			r(X_s,A_s,B_s)+\alpha(X_s)J_{\delta_\varepsilon}(X_{s+1})
			\mid\mathcal G_s\right]
			+\bar\alpha\delta_\varepsilon. \tag{5.15}
		\end{align}
		Fix a deterministic $T\in\mathbb N$. For $0\leq s\leq T-1$, the product $\left(\prod_{k=0}^{s-1}\alpha(X_k)\right)I_{\{s<T\wedge\theta_\varepsilon^1\wedge\theta^2\wedge\sigma_D\}}$
		is $\mathcal G_s$-measurable. Hence,
		\begin{align*}
			&E_x^{\pi_\varepsilon^1,\pi^2}\!\bigg[
			\left(\prod_{k=0}^{s-1}\alpha(X_k)\right)I_{\{s<T\wedge\theta_\varepsilon^1\wedge\theta^2\wedge\sigma_D\}}
			E_x^{\pi_\varepsilon^1,\pi^2}\!\left[
			r(X_s,A_s,B_s)+\alpha(X_s)J_{\delta_\varepsilon}(X_{s+1})
			\mid\mathcal G_s\right]\bigg]\\
			&=
			E_x^{\pi_\varepsilon^1,\pi^2}\!\bigg[
			\left(\prod_{k=0}^{s-1}\alpha(X_k)\right)I_{\{s<T\wedge\theta_\varepsilon^1\wedge\theta^2\wedge\sigma_D\}}
			\left(r(X_s,A_s,B_s)+\alpha(X_s)J_{\delta_\varepsilon}(X_{s+1})\right)\bigg].
		\end{align*}
		Multiply (5.15) by $\left(\prod_{k=0}^{s-1}\alpha(X_k)\right) I_{\{s<T\wedge\theta_\varepsilon^1\wedge\theta^2\wedge\sigma_D\}}$, take expectations, and sum from $s=0$ to $T-1$. Using
		\begin{align*}
			&\sum_{s=0}^{T\wedge\theta_\varepsilon^1\wedge\theta^2\wedge\sigma_D-1}\left(\prod_{k=0}^{s-1}\alpha(X_k)\right)J_{\delta_\varepsilon}(X_s)
			-\sum_{s=0}^{T\wedge\theta_\varepsilon^1\wedge\theta^2\wedge\sigma_D-1}\left(\prod_{k=0}^{s}\alpha(X_k)\right)J_{\delta_\varepsilon}(X_{s+1})\\
			&\qquad=J_{\delta_\varepsilon}(x)
			-\left(\prod_{k=0}^{T\wedge\theta_\varepsilon^1\wedge\theta^2\wedge\sigma_D-1}\alpha(X_k)\right)
			J_{\delta_\varepsilon}(X_{T\wedge\theta_\varepsilon^1\wedge\theta^2\wedge\sigma_D}),
		\end{align*}
		and the bound
		\[
			\bar\alpha\delta_\varepsilon
			E_x^{\pi_\varepsilon^1,\pi^2}\!\left[
			\sum_{s=0}^{T\wedge\theta_\varepsilon^1\wedge\theta^2\wedge\sigma_D-1}
			\left(\prod_{k=0}^{s-1}\alpha(X_k)\right)\right]
			\leq\bar\alpha\delta_\varepsilon\sum_{s=0}^{\infty}\bar\alpha^s
			=\frac{\bar\alpha\delta_\varepsilon}{1-\bar\alpha},
		\]
		we obtain
		\begin{align}
			J_{\delta_\varepsilon}(x)
			&\leq E_x^{\pi_\varepsilon^1,\pi^2}\bigg[
			\sum_{s=0}^{T\wedge\theta_\varepsilon^1\wedge\theta^2\wedge\sigma_D-1}\left(\prod_{k=0}^{s-1}\alpha(X_k)\right)
			r(X_s,A_s,B_s)\bigg]\notag\\
			&\quad+E_x^{\pi_\varepsilon^1,\pi^2}\bigg[
			\left(\prod_{k=0}^{T\wedge\theta_\varepsilon^1\wedge\theta^2\wedge\sigma_D-1}\alpha(X_k)\right)
			J_{\delta_\varepsilon}(X_{T\wedge\theta_\varepsilon^1\wedge\theta^2\wedge\sigma_D})\bigg]
			+\frac{\bar\alpha\delta_\varepsilon}{1-\bar\alpha}. \tag{5.16}\label{eq:epsilon-truncated-1}
		\end{align}

		On $\{\theta_\varepsilon^1<\sigma_D\}$, (5.10)--(5.11) give
		\begin{equation}
			J_{\delta_\varepsilon}(X_{\theta_\varepsilon^1})
			\leq g_1(X_{\theta_\varepsilon^1})+\frac{\varepsilon}{4}.\tag{5.17}\label{eq:epsilon-stop-1}
		\end{equation}
		On $\{\theta^2<\theta_\varepsilon^1,\ \theta^2<\sigma_D\}$, one has $X_{\theta^2}\notin D_\varepsilon^1$, and hence
		\[
			g_1(X_{\theta^2})
			<J_{\delta_\varepsilon}(X_{\theta^2})-\frac{\varepsilon}{4}
			<J(X_{\theta^2})
			\leq\max\{g_1(X_{\theta^2}),g_2(X_{\theta^2})\}.
		\]
		Consequently, $g_2(X_{\theta^2})>g_1(X_{\theta^2})$, so
		\begin{equation}
			J_{\delta_\varepsilon}(X_{\theta^2})
			\leq J(X_{\theta^2})
			\leq g_2(X_{\theta^2}).\tag{5.18}\label{eq:epsilon-opponent-stop-2}
		\end{equation}
		For every $\theta\in\Theta$, the discount factors also give
		\begin{equation}
			0\leq\left(\prod_{k=0}^{\theta-1}\alpha(X_k)\right)
			I_{\{\theta<\sigma_D\}}\leq1.\tag{5.19}\label{eq:epsilon-discount-bound}
		\end{equation}

		By the same argument as in the proof of Theorem 4.2, using the dominated convergence theorem with the corresponding function and strategies, we may let $T\to\infty$ in (5.16). The terminal term is zero on $\{\sigma_D\leq\theta_\varepsilon^1\wedge\theta^2,\ \sigma_D<\infty\}$ because $J_{\delta_\varepsilon}=0$ on $D$. On the remaining events, splitting according to which player stops first and using the inequalities above gives
		\begin{align}
			J_{\delta_\varepsilon}&(x)
			\leq E_x^{\pi_\varepsilon^1,\pi^2}\bigg[
			\sum_{s=0}^{\theta_\varepsilon^1\wedge\theta^2\wedge\sigma_D-1}\left(\prod_{k=0}^{s-1}\alpha(X_k)\right)
			r(X_s,A_s,B_s)+\left(\prod_{k=0}^{\theta_\varepsilon^1-1}\alpha(X_k)\right)
			J_{\delta_\varepsilon}(X_{\theta_\varepsilon^1})
			I_{\{\theta_\varepsilon^1\leq\theta^2,\,\theta_\varepsilon^1<\sigma_D\}}\notag\\
			&\hspace{40mm}+\left(\prod_{k=0}^{\theta^2-1}\alpha(X_k)\right)J_{\delta_\varepsilon}(X_{\theta^2})
			I_{\{\theta^2<\theta_\varepsilon^1,\,\theta^2<\sigma_D\}}
			\bigg]+\frac{\bar\alpha\delta_\varepsilon}{1-\bar\alpha}\quad\text{by }\eqref{eq:epsilon-truncated-1}\notag\displaybreak[1]\\
			&\leq E_x^{\pi_\varepsilon^1,\pi^2}\bigg[
			\sum_{s=0}^{\theta_\varepsilon^1\wedge\theta^2\wedge\sigma_D-1}\left(\prod_{k=0}^{s-1}\alpha(X_k)\right)
			r(X_s,A_s,B_s)+\left(\prod_{k=0}^{\theta_\varepsilon^1-1}\alpha(X_k)\right)
			\left(g_1(X_{\theta_\varepsilon^1})+\frac{\varepsilon}{4}\right)
			I_{\{\theta_\varepsilon^1\leq\theta^2,\,\theta_\varepsilon^1<\sigma_D\}}\notag\\
			&\hspace{40mm}+\left(\prod_{k=0}^{\theta^2-1}\alpha(X_k)\right)g_2(X_{\theta^2})
			I_{\{\theta^2<\theta_\varepsilon^1,\,\theta^2<\sigma_D\}}
			\bigg]+\frac{\bar\alpha\delta_\varepsilon}{1-\bar\alpha}\quad\text{by }\eqref{eq:epsilon-stop-1},\eqref{eq:epsilon-opponent-stop-2}\notag\displaybreak[1]\\
			&=J(x,\pi_\varepsilon^1,\theta_\varepsilon^1,\pi^2,\theta^2)\
			+\frac{\varepsilon}{4}
			E_x^{\pi_\varepsilon^1,\pi^2}\!\left[
			\left(\prod_{k=0}^{\theta_\varepsilon^1-1}\alpha(X_k)\right)
			I_{\{\theta_\varepsilon^1\leq\theta^2,\,\theta_\varepsilon^1<\sigma_D\}}\right]
			+\frac{\bar\alpha\delta_\varepsilon}{1-\bar\alpha}\quad\text{by (2.6)}\notag\displaybreak[1]\\
			&\leq J(x,\pi_\varepsilon^1,\theta_\varepsilon^1,\pi^2,\theta^2)
			+\frac{\varepsilon}{4}+\frac{\bar\alpha\delta_\varepsilon}{1-\bar\alpha}\quad\text{by }\eqref{eq:epsilon-discount-bound}.\tag{5.20}\label{eq:epsilon-criterion-1}
		\end{align}
		Using $J\leq J_{\delta_\varepsilon}+\delta_\varepsilon$ and $\frac{\delta_\varepsilon}{(1-\bar\alpha)}=\frac{\bar\alpha\varepsilon}{6}$, we obtain
		\begin{align}
			J(x,\pi_\varepsilon^1,\theta_\varepsilon^1,\pi^2,\theta^2)
			&\geq J(x)-\left(\frac14+\frac{\bar\alpha}{6}\right)\varepsilon
			\geq J(x)-\varepsilon. \tag{5.21}
		\end{align}

		On the other hand, fix $(\pi^1,\theta^1)\in\Pi^1\times\Theta$. On $\{s<\theta_\varepsilon^2\wedge\sigma_D\}$, (5.8), (5.10), and $\delta_\varepsilon<\frac{\varepsilon}{4}$ give $J(X_s)<g_2(X_s)$. As above, Lemma 4.1(b), (5.12), and (5.14) yield
		\begin{equation}
			J_{\delta_\varepsilon}(X_s)
			\geq E_x^{\pi^1,\pi_\varepsilon^2}\!\left[
			r(X_s,A_s,B_s)+\alpha(X_s)J_{\delta_\varepsilon}(X_{s+1})
			\mid\mathcal G_s\right]-\delta_\varepsilon,
			\quad P_x^{\pi^1,\pi_\varepsilon^2}\text{-a.s.}\tag{5.22}\label{eq:epsilon-step-2}
		\end{equation}
		By (5.9) and (5.10)--(5.11), respectively,
		\begin{align}
			J_{\delta_\varepsilon}(X_{\theta^1})&\geq g_1(X_{\theta^1})
			\quad\text{on }\{\theta^1<\sigma_D\},\tag{5.23}\label{eq:epsilon-opponent-stop-1}\\
			J_{\delta_\varepsilon}(X_{\theta_\varepsilon^2})&\geq g_2(X_{\theta_\varepsilon^2})-\frac{\varepsilon}{4}
			\quad\text{on }\{\theta_\varepsilon^2<\sigma_D\}.\tag{5.24}\label{eq:epsilon-stop-2}
		\end{align}
		Proceeding as above, we sum the discounted inequalities in \eqref{eq:epsilon-step-2} and let $T\to\infty$. The accumulated error is at most $\delta_\varepsilon\sum_{s=0}^\infty\bar\alpha^s=\frac{\delta_\varepsilon}{1-\bar\alpha}$. Using \eqref{eq:epsilon-opponent-stop-1}--\eqref{eq:epsilon-stop-2}, \eqref{eq:epsilon-discount-bound}, and the zero extension on $D$ gives
		\begin{align}
			J_{\delta_\varepsilon}(x)
			&\geq J(x,\pi^1,\theta^1,\pi_\varepsilon^2,\theta_\varepsilon^2)
			-\frac{\varepsilon}{4}
			E_x^{\pi^1,\pi_\varepsilon^2}\!\left[
			\left(\prod_{k=0}^{\theta_\varepsilon^2-1}\alpha(X_k)\right)
			I_{\{\theta_\varepsilon^2<\theta^1,\,\theta_\varepsilon^2<\sigma_D\}}\right]
			-\frac{\delta_\varepsilon}{1-\bar\alpha}\notag\\
			&\geq J(x,\pi^1,\theta^1,\pi_\varepsilon^2,\theta_\varepsilon^2)
			-\frac{\varepsilon}{4}-\frac{\delta_\varepsilon}{1-\bar\alpha}
			\quad\text{by }\eqref{eq:epsilon-discount-bound}.\tag{5.25}\label{eq:epsilon-criterion-2}
		\end{align}
		Hence,
		\begin{align}
			J(x,\pi^1,\theta^1,\pi_\varepsilon^2,\theta_\varepsilon^2)
			\leq J_{\delta_\varepsilon}(x)
			+\frac{\varepsilon}{4}+\frac{\delta_\varepsilon}{1-\bar\alpha}
			\leq J(x)+\left(\frac{1}{4}+\frac{\bar\alpha}{6}\right)\varepsilon
			\leq J(x)+\varepsilon. \tag{5.26}
		\end{align}
		By (5.21) and (5.26), for every $x\in D^c$,
		\begin{align*}
			\inf_{(\pi^2,\theta^2)\in\Pi^2\times\Theta}
			J(x,\pi_\varepsilon^1,\theta_\varepsilon^1,\pi^2,\theta^2)
			\geq J(x)
			-\left(\frac{1}{4}+\frac{\bar\alpha}{6}\right)\varepsilon
			\geq J(x)-\varepsilon,
		\end{align*}
		and
		\begin{align*}
			\sup_{(\pi^1,\theta^1)\in\Pi^1\times\Theta}
			J(x,\pi^1,\theta^1,\pi_\varepsilon^2,\theta_\varepsilon^2)
			\leq J(x)
			+\left(\frac{1}{4}+\frac{\bar\alpha}{6}\right)\varepsilon
			\leq J(x)+\varepsilon.
		\end{align*}
		Thus, by (2.11)--(2.12),
		$(\pi_\varepsilon^1,\theta_\varepsilon^1;
		\pi_\varepsilon^2,\theta_\varepsilon^2)$
		is an $\varepsilon$-Nash equilibrium.
	\end{proof}
	The proof of Theorem 5.3 leads to the following algorithm.
	
	\begin{samepage}
	\textbf{Algorithm \RN{2} (Computation of an $\varepsilon$-Nash equilibrium):}\\
	Step 1 (Initialization): Input an accuracy $\varepsilon>0$, set $\delta_\varepsilon:=\frac{\bar\alpha(1-\bar\alpha)\varepsilon}{6}$, and run Algorithm \RN{1} with accuracy $\delta_\varepsilon$ to obtain $J_{\delta_\varepsilon}$ satisfying (5.8).\\
	Step 2 (Stopping sets): Construct $D_\varepsilon^1:=\{x\in D^c:g_1(x)-J_{\delta_\varepsilon}(x)\geq -\frac{\varepsilon}{4}\}$ and $D_\varepsilon^2:=\{x\in D^c:g_2(x)-J_{\delta_\varepsilon}(x)\leq \frac{\varepsilon}{4}\}$.\\
	Step 3 ($\varepsilon$-optimal policies): Choose arbitrary admissible stationary controls on $D$. For each $x\in D^c$, solve (5.3) and (5.4) with $u(x):=J_{\delta_\varepsilon}(x)$ to obtain $\pi_\varepsilon^1(\cdot\mid x)$ and $\pi_\varepsilon^2(\cdot\mid x)$.\\
	Step 4 ($\varepsilon$-optimal stopping): Define $\theta_\varepsilon^1:=\inf\{t\geq0:X_t\in D_\varepsilon^1\}$ and $\theta_\varepsilon^2:=\inf\{t\geq0:X_t\in D_\varepsilon^2\}$. System failure terminates the game at $\sigma_D$ regardless of these stopping decisions. Then $(\pi_\varepsilon^1,\theta_\varepsilon^1;\pi_\varepsilon^2,\theta_\varepsilon^2)$ is an $\varepsilon$-Nash equilibrium.
	\end{samepage}
	\begin{corollary}
		If $X$, $A(x)$, and $B(x)$ are finite, then, for every $\varepsilon>0$, Algorithms \RN{1} and \RN{2} produce an $\varepsilon$-Nash equilibrium $(\pi_\varepsilon^1,\theta_\varepsilon^1;\pi_\varepsilon^2,\theta_\varepsilon^2)$.
	\end{corollary}
	
	\begin{remark}
		Ghosh and Mallikarjuna Rao \autocite{Ghosh_2007} characterize exact equilibria through the value equation but do not provide a numerical procedure with a prescribed equilibrium accuracy. In the present model, Theorem 5.2 gives, for finite $X_*$ and $n\geq1$,
		\[
			0\leq J(x)-\hat u_n(x)
			\leq\frac{\bar\alpha}{1-\bar\alpha}
			\max_{y\in X_*}|\hat u_n(y)-\hat u_{n-1}(y)|.
		\]
		The contraction-based value-iteration principle is classical. Here its error bound is linked, through the choice of $\delta_\varepsilon$, to both players' stationary controls and stopping times with an $\varepsilon$-Nash guarantee. For finite state and action spaces, the procedure uses linear programs with the state-dependent discount factors and updates only one value sequence.
	\end{remark}

	\section{A queueing game with unbounded rewards}\label{s:6}

This section considers a service queue. We first specify the model and verify the conditions for the existence of the value and a Nash equilibrium. We then choose numerical parameters to illustrate the computation of the value and an $\varepsilon$-Nash equilibrium.

Consider a client who submits tasks to a service provider under an agreement with service charges, compensation for unfinished tasks, and termination payments. The client is player 1, and the provider is player 2. There is no upper bound on the number of unfinished tasks, including any task being processed. If $x$ is the number of unfinished tasks, the state space is $X:=\{0,1,2,\ldots\}$. A positive reward is a payment from the provider to the client, and a negative reward is a payment in the opposite direction. Thus, the client maximizes the expected discounted criterion in (2.6), whereas the provider minimizes the same criterion. Either party may terminate the agreement. We take $D=\varnothing$, so termination is determined by their stopping decisions rather than system failure.

At each decision epoch when both parties continue, the client chooses whether to increase its submission frequency: $a=0$ means normal submission, and $a=1$ means increased submission. The provider chooses whether to allocate additional processing capacity: $b=0$ means normal processing, and $b=1$ means additional processing capacity. Thus, $A(x)=B(x):=\{0,1\}$ for every $x\in X$. We assume that at most one task is submitted and at most one task is completed during each period, and that submission and completion do not occur in the same period. Given actions $a$ and $b$, a new task is submitted and admitted with probability $q^+(a)$.
If $x>0$, one task is completed with probability $q^-(b)$.
If $x=0$, no task can be completed.
Assume that
\[
0\leq q^+(0)<q^+(1),\qquad
0\leq q^-(0)<q^-(1),\qquad
q^+(1)+q^-(1)\leq1.
\]
Thus, for $x\geq1$, the queue length increases by one with probability
$q^+(a)$, decreases by one with probability $q^-(b)$, and remains
unchanged with probability $1-q^+(a)-q^-(b)$. Consequently, the transition probabilities are
\begin{align}
	p(y\mid x,a,b)
	={}&q^+(a)I_{\{y=x+1\}}
	+q^-(b)I_{\{y=\max\{x-1,0\}\}}
	+\bigl(1-q^+(a)-q^-(b)\bigr)I_{\{y=x\}}. \tag{6.1}
\end{align}
At the empty state, (6.1) gives
\[
	p(1\mid0,a,b)=q^+(a),\qquad
	p(0\mid0,a,b)=1-q^+(a).
\]

Let $c,c_a,c_b$ be nonnegative constants and $h>0$. In each period of continued operation, the client pays a basic service charge $c$ and receives compensation $h$ per unfinished task. The client also pays an additional submission fee $c_a$ if $a=1$ and an additional service fee $c_b$ if $b=1$. These are payments between the two parties; external operating costs are not included in the criterion. The resulting one-stage reward to the client is
\begin{equation}
	r(x,a,b):=-c+hx-c_a a-c_b b. \tag{6.2}
\end{equation}
Let $\eta,\zeta,k_1$ be nonnegative constants and $k_2>0$, describing the termination payments. The agreement specifies an integer $N\geq1$: when fewer than $N$ tasks remain and the client terminates the agreement, a fixed charge $\eta$ is reduced by a refund $k_1$ per unfinished task. If the provider terminates it first, the client receives compensation $\zeta$ plus $k_2$ per unfinished task. When at least $N$ tasks remain, the client receives this latter compensation regardless of which party terminates the agreement. Thus $N$ specifies when the common termination payment applies; it is not a capacity limit. The terminal payoff functions are
\begin{equation}
	g_1(x):=\begin{cases}
	-\eta+k_1x,&0\leq x<N,\\
	\zeta+k_2x,&x\geq N,
	\end{cases}
	\qquad
	g_2(x):=\zeta+k_2x. \tag{6.3}
\end{equation}
As in (2.6), $g_1$ applies when player 1 stops no later than player 2, and $g_2$ applies when player 2 stops first. All rewards are measured from the client's perspective. No ordering between $g_1$ and $g_2$ is imposed for $x<N$. For $x\geq N$, the two termination payments agree. Since $h,k_2>0$, the running reward and both terminal rewards tend to $+\infty$ as $x\to\infty$.

Both parties evaluate future net payments using a common discount factor $\alpha(x)\in(0,1)$ at state $x$. Assume that $\bar\alpha=\sup_{x\in X}\alpha(x)<1$. The cumulative discount is the product of these factors along the state path, as in (2.6).

\begin{proposition}
	Under the conditions specified above, the model in (6.1)--(6.3) satisfies Assumptions 3.1 and 3.3 with $V_0(x)=\ell+x$ for a sufficiently large constant $\ell\geq1$. The value and a Nash equilibrium exist, and Algorithms \RN{1} and \RN{2} apply.
\end{proposition}
\begin{proof}
	The conditions on $q^+$ and $q^-$ give $0\leq q^+(a)+q^-(b)\leq1$ for every $(a,b)$. Hence (6.1) defines a transition probability distribution on $X$. Since $\bar\alpha<1$, choose $\ell\geq1$ such that $\bar\alpha\left(1+\frac{q^+(1)}{\ell}\right)<1.$
	Taking $V_0(x)=\ell+x$, we obtain
	\begin{align*}
		\sum_{y\in X}V_0(y)p(y\mid x,a,b)
		=V_0(x)+q^+(a)-q^-(b)I_{\{x\geq1\}}
		\leq V_0(x)+q^+(1)
		\leq\left(1+\frac{q^+(1)}{\ell}\right)V_0(x).
	\end{align*}
	This proves (3.1) with $m=1+\frac{q^+(1)}{\ell}$ and $\bar\alpha m<1$. Set
	\[
		c_0:=\max\left\{\frac{c+c_a+c_b}{\ell},h,
		\frac{\eta}{\ell},k_1,\frac{\zeta}{\ell},k_2\right\}.
	\]
	By (6.2)--(6.3), for every $(x,a,b)$,
	\begin{align*}
		|r(x,a,b)|&\leq c+c_a+c_b+hx\leq c_0V_0(x),\\
		|g_1(x)|&\leq\max\{\eta+k_1x,\zeta+k_2x\}\leq c_0V_0(x),\\
		|g_2(x)|&=\zeta+k_2x\leq c_0V_0(x).
	\end{align*}
	Thus Assumption 3.1 holds. Finiteness of the action sets gives compactness and the continuity required in Assumption 3.3. Moreover,
	\[
	X_*:=\{x\in X:g_1(x)\neq g_2(x)\}\subseteq\{0,1,\ldots,N-1\},
	\]
	so $X_*$ is finite. Theorem 4.2 and Theorems 5.2--5.3 give the stated conclusions. In particular,
	\[
	\hat u_n(x)=J(x)=\zeta+k_2x,\qquad x\geq N,\quad n\geq0.
	\]
\end{proof}

\begin{remark}
	In this model, $r$, $g_1$, and $g_2$ are unbounded, unlike the bounded-reward assumptions in the main results of Ghosh and Mallikarjuna Rao \autocite{Ghosh_2007}, which also use a constant discount factor. Pal and Saha \autocite{Pal_Saha_2020} give a continuous-time queueing example with a constant discount rate. Here, the decisions are made in discrete time and the discount factor depends on the current queue length. The common termination payment for $x\geq N$ makes $X_*$ finite without bounding the rewards.
\end{remark}

For the numerical calculation, we take $N=7,\quad q^+(a)=0.22+0.16a,\quad q^-(b)=0.22+0.16b,\quad
c=1.15,\quad h=0.38,\quad c_a=c_b=0.22,\quad
\eta=1.55,\quad k_1=0.22,\quad \zeta=0.35,\quad k_2=0.48,
\alpha(x)=0.80+0.12\frac{x}{1+x}.$ Here $N=7$ means that the common termination payment applies when at least seven tasks remain; tasks may still be submitted at these states if both parties continue. The discount factors increase with the number of unfinished tasks and remain below $0.92$. These numerical parameters illustrate the agreement and are not fitted to operational data. In particular, $q^+(1)+q^-(1)=0.76<1$ and $\bar\alpha=\sup_{x\in X}\alpha(x)=0.92$.
For these parameters, $g_1(x)<g_2(x)$ for $0\leq x\leq6$, whereas $g_1(x)=g_2(x)$ for $x\geq7$. With $\ell=10$ in Proposition 6.1, $V_0(x)=10+x$, $m=1.038$, $c_0=0.48$, $\bar\alpha m=0.95496<1$.

Although the running and terminal rewards are unbounded, their absolute values are bounded by $0.48V_0(x)$. Also, $X_*=\{0,1,\ldots,6\}$, $\hat u_n(x)=J(x)=0.35+0.48x\quad(x\geq7,\ n\geq0).$

Applying Algorithm \RN{1} with $\hat u_0=g_1$ and $\delta=10^{-6}$, the stopping rule is first satisfied at $n=49$. With $J_\delta:=\hat u_{49}$, Theorem 5.2 gives
\[
0\leq J(x)-J_\delta(x)
\leq\frac{\bar\alpha}{1-\bar\alpha}
\max_{y\in X_*}|\hat u_{49}(y)-\hat u_{48}(y)|
\approx7.778\times10^{-7}<\delta,\qquad x\in X.
\]
Table~\ref{tab:project-values} gives $J_\delta$ on $X_*$, rounded to six decimal places, together with $g_1$ and $g_2$. For all $x\geq7$, $J_\delta(x)=g_1(x)=g_2(x)=0.35+0.48x$ exactly.

\begin{table}[htbp]
	\centering
	\caption{Approximation $J_\delta$ on $X_*$ for $\delta=10^{-6}$ after 49 iterations.}
	\label{tab:project-values}
	\begin{tabular}{crrrrrrr}
		\toprule
		State $x$ & 0 & 1 & 2 & 3 & 4 & 5 & 6\\
		\midrule
		$g_1(x)$
		& $-1.550000$ & $-1.330000$ & $-1.110000$ & $-0.890000$ & $-0.670000$ & $-0.450000$ & $-0.230000$\\
		$J_\delta(x)$
		& $-1.550000$ & $-1.330000$ & $-1.110000$ & $-0.798776$ & $0.464711$ & $2.075736$ & $3.230000$\\
		$g_2(x)$
		& $0.350000$ & $0.830000$ & $1.310000$ & $1.790000$ & $2.270000$ & $2.750000$ & $3.230000$\\
		\bottomrule
	\end{tabular}
\end{table}

We next apply Algorithm \RN{2} with $\varepsilon=10^{-4}$. Its required value-approximation accuracy is
\[
	\delta_\varepsilon=\frac{\bar\alpha(1-\bar\alpha)}{6}\varepsilon
	=\frac{92}{75}\times10^{-6}>\delta.
\]
Thus the approximation already computed satisfies (5.8), and we may take $J_{\delta_\varepsilon}:=J_\delta$. The stopping sets defined by (5.10) and the states where both players continue are
\begin{align}
	D_\varepsilon^1&=\{0,1,2\}\cup\{x\in X:x\geq7\},\qquad
	D_\varepsilon^2=\{x\in X:x\geq6\},\notag\\
	C&:=X\setminus(D_\varepsilon^1\cup D_\varepsilon^2)=\{3,4,5\}. \tag{6.4}
\end{align}
For these parameters, the client terminates the agreement when at most two or at least seven tasks remain unfinished. The provider terminates it when at least six tasks remain. Both parties continue when there are three, four, or five unfinished tasks. Solving (5.3)--(5.4) with $u=J_{\delta_\varepsilon}$ gives the stationary policies $\pi_\varepsilon^1$ and $\pi_\varepsilon^2$. Table~\ref{tab:project-policies} gives their action probabilities on $C$. In this example, these probabilities are zero or one: the client increases submission only at $x=4$, and the provider allocates additional processing capacity at every state in $C$.

\begin{table}[htbp]
	\centering
	\caption{Stationary policies from Algorithm \RN{2} on $C$ for $\varepsilon=10^{-4}$.}
	\label{tab:project-policies}
	\begin{tabular}{ccccc}
		\toprule
		$x$ & $\pi^1_\varepsilon(0\mid x)$ & $\pi^1_\varepsilon(1\mid x)$
		& $\pi^2_\varepsilon(0\mid x)$ & $\pi^2_\varepsilon(1\mid x)$\\
		\midrule
		3 & 1 & 0 & 0 & 1\\
		4 & 0 & 1 & 0 & 1\\
		5 & 1 & 0 & 0 & 1\\
		\bottomrule
	\end{tabular}
\end{table}

For $x\geq7$, the known values give $J_{\delta_\varepsilon}(x+1)-J_{\delta_\varepsilon}(x)=0.48$, $J_{\delta_\varepsilon}(x-1)-J_{\delta_\varepsilon}(x)\leq-0.48.$
Consequently, for $a,b\in\{0,1\}$,
\begin{align*}
\mathbb I_\alpha J_{\delta_\varepsilon}(x,1,b)
-\mathbb I_\alpha J_{\delta_\varepsilon}(x,0,b)
&=-0.22+0.0768\alpha(x)<0,\\
\mathbb I_\alpha J_{\delta_\varepsilon}(x,a,1)
-\mathbb I_\alpha J_{\delta_\varepsilon}(x,a,0)
&\leq-0.22-0.0768\alpha(x)<0.
\end{align*}
Thus the solutions of (5.3)--(5.4) at every $x\geq7$ select $a=0$ and $b=1$, without requiring a separate numerical calculation at each such state. The finite calculations at $x\in\{0,1,2,6\}$ give the same controls. Hence the policies are specified on all of $X$ by
\[
\pi_\varepsilon^1(1\mid x)=I_{\{x=4\}},\qquad
\pi_\varepsilon^2(1\mid x)=1,\qquad x\in X.
\]
Both players stop immediately when $x\geq7$, so these controls also define the state process if it is continued after stopping.

The approximation satisfies (5.8), and the policies satisfy (5.12). For the stopping sets in (6.4), (5.11) gives
\begin{align*}
	\theta_\varepsilon^1&=\inf\{t\geq0:X_t\leq2\ \text{or}\ X_t\geq7\},\\
	\theta_\varepsilon^2&=\inf\{t\geq0:X_t\geq6\}.
\end{align*}
Consequently, Theorem 5.3 shows that $(\pi_\varepsilon^1,\theta_\varepsilon^1;\pi_\varepsilon^2,\theta_\varepsilon^2)$ is a $10^{-4}$-Nash equilibrium.

These stopping times depend on the state path. For illustration, suppose that initially four tasks are unfinished, so that $X_0=4$, and consider the following two possible paths under $(\pi_\varepsilon^1,\pi_\varepsilon^2)$:
\[
	\begin{array}{c|rrrrrrr|cc}
		t & 0 & 1 & 2 & 3 & 4 & 5 & 6
		& \theta_\varepsilon^1 & \theta_\varepsilon^2\\ \hline
		X_t & 4 & 3 & 2 & 3 & 4 & 5 & 6 & 2 & 6\\
		X_t & 4 & 5 & 6 & 5 & 4 & 3 & 2 & 6 & 2
	\end{array}
\]
Each displayed transition has positive probability by (6.1). The game ends at time $2$ on both paths; the later entries describe the state process continued under the same policies, not further play. The stopping-time values are path dependent.

\begin{figure}[H]
	\centering
	\includegraphics[width=\linewidth]{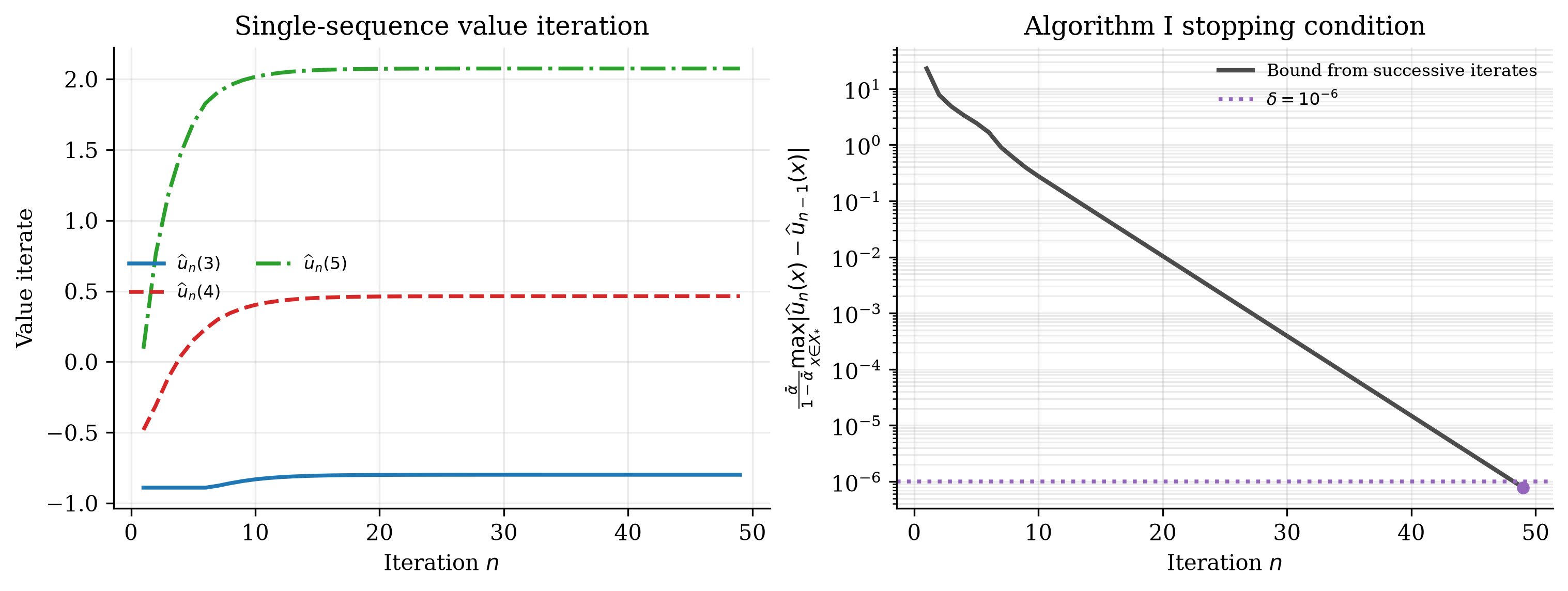}
	\caption{The sequence $\hat u_n$ (left) and the value-error bound from successive iterates (right).}
	\label{fig:project-convergence}
\end{figure}

Figure~\ref{fig:project-convergence} illustrates the single-sequence convergence in Theorem 5.1 and the stopping condition in Algorithm \RN{1}. The left panel shows $\hat u_n(x)$ for $x\in C$. The right panel shows $\frac{\bar\alpha}{1-\bar\alpha}\max_{x\in X_*}|\hat u_n(x)-\hat u_{n-1}(x)|$ on a logarithmic scale; the horizontal line marks $\delta=10^{-6}$.

\end{document}